\documentclass[UTF8,a4paper,12pt]{article}

\usepackage[T1]{fontenc} 
\usepackage{hyperref}
\hypersetup{
	colorlinks=true,
	linkcolor=blue,
	filecolor=blue,      
	urlcolor=black,
	citecolor=blue,
	anchorcolor=black,
	menucolor=red
}
\usepackage{dsfont}
\usepackage{amsmath}
\usepackage{amsthm}
\usepackage{amssymb}
\usepackage{mathrsfs}
\usepackage[english]{babel}
\usepackage{comment}
\usepackage{enumerate}
\usepackage{geometry}
\usepackage{pythonhighlight}
\usepackage{graphicx}
\usepackage{fancyhdr}
\usepackage{authblk}
\usepackage{url}
\usepackage{float}
\usepackage{booktabs}
\usepackage[all]{xy}
\usepackage{color}
\usepackage{natbib}
\usepackage{hyperref}

\theoremstyle{plain}
\newtheorem{theorem}{Theorem}

\newtheorem{lemma}{Lemma}
\newtheorem{proposition}{Proposition}

\theoremstyle{definition}
\newtheorem{definition}{Definition}

\newtheorem{assumption}{Assumption}

\theoremstyle{remark}
\newtheorem{remark}{Remark}
\newtheorem{problem}{Problem}

\newcommand{\ud}{\mathrm{d}}

\newcommand{\E}{\mathbb{E}}

\newcommand{\F}{\mathcal{F}}
\newcommand{\e}{\mathcal{E}}

\newcommand{\R}{\mathbb{R}}

\newcommand{\p}{\mathbb{P}}
\newcommand{\Q}{\mathbb{Q}}

\newcommand{\one}{\mathbf{1}}

\renewcommand{\(}{\left(}
\renewcommand{\)}{\right)}
\renewcommand{\[}{\left[}
\renewcommand{\]}{\right]}

\renewcommand{\geq}{\geqslant}
\renewcommand{\leq}{\leqslant}
\renewcommand{\epsilon}{\varepsilon}

\newcommand\keywords[1]{\textbf{Keywords}: #1}

\title{Comparison Between Mean-Variance and Monotone Mean-Variance Preferences Under Jump Diffusion and Stochastic Factor Model}

\author[a]{Yuchen Li}
\author[a]{Zongxia Liang}
\author[b]{Shunzhi Pang}
\affil[a]{\small{Department of Mathematical Sciences, Tsinghua University, Beijing 100084, China}}
\affil[b]{\small{School of Economics and Management, Tsinghua University, Beijing 100084, China}}

\date{}
\begin{document}
	\maketitle
	
\def\itNo#1#2#3{\int_{\mathbb{R}}#1 \tilde{N}\(\ud #2,\ud #3\)}
\def\itNt#1#2#3{\int_{\mathbb{R}}#1 \tilde{N}^2\(\ud #2,\ud #3\)}
\begin{abstract}
	This paper compares the optimal investment problems based on monotone mean-variance (MMV) and mean-variance (MV) preferences in the L\'{e}vy market with an untradable stochastic factor. It is an open question proposed by Trybu{\l}a and Zawisza. Using the dynamic programming and Lagrange multiplier methods, we get the HJBI and HJB equations corresponding to the two investment problems. The equations are transformed into a new-type parabolic equation, from which the optimal strategies under both preferences are derived. We prove that the two optimal strategies and value functions coincide if and only if an important market assumption holds. When the assumption violates, MMV investors act differently from MV investors. Thus, we conclude that the difference between continuous-time MMV and MV portfolio selections is due to the discontinuity of the market. In addition, we derive the efficient frontier and analyze the economic impact of the jump diffusion risky asset. We also provide empirical evidences to demonstrate the validity of the assumption in real financial market.
\end{abstract}

\section{Introduction.}

Since \citet{markowits1952portfolio} proposed the mean-variance preference, it became a popular criterion in portfolio selection problems. Previous works have widely considered both multi-period discrete and continuous-time portfolio selection problems under mean-variance preferences. \citet{LN2000optimal} and \citet{ZL2000continuous} use embedding techniques to overcome the time-inconsistency of mean-variance preferences, leading a great number of following works, see e.g., \citet{cerny2007structure}, \citet{zhou2002noshort}, \citet{lim2002random} and \citet{lim2004random}. However, when we consider the optimal investment problems, there is still another drawback of mean-variance preferences: non-monotonicity.

For a preference $\rho$ with monotonicity, if two random variables satisfy $X\leq Y$ almost surely, we expect $\rho(X)\leq \rho(Y)$. However, mean-variance preferences do not have such property. A strictly wealthier investor may have a lower utility than others, which violates the basic assumption of economic rationality. \citet{maccheroni2009portfolio} provide an example to illustrate the non-monotonicity of mean-variance preferences. To conquer such a disadvantage, they propose a new criterion named the monotone mean-variance preference. It is the best possible monotone approximation of the classical mean-variance preference without loss of good tractability. In this paper, we use the abbreviations MMV and MV to denote monotone mean-variances and mean-variances. 

There are several works related to MMV preferences. \citet{trybula2019continuous} study the optimal investment problem in a continuous-time financial market model with a stochastic factor. Using the method of dynamic programming and solving the Hamilton–Jacobi–Bellman–Isaacs (HJBI) equation, they propose a verification theorem and obtain an explicit optimal strategy of the investment under MMV preferences. \citet{trybula2019continuous} use the result from \citet{zawisza2012quadratic} and prove that the optimal strategies under MMV and MV preferences are the same, which is different from the conclusion of the single-period case considered in \citet{maccheroni2009portfolio}. The reason of this difference may be either the discreteness or the discontinuity of the wealth process. Using the analytic approach rather than comparing the two optimal strategies directly, \citet{SL2020note} generalize the conclusion of \citet{trybula2019continuous} to all continuous semimartingale financial market models. In the continuous market scenario, \citet{LG2021RAIRO-OR} solve the optimal reinsurance and investment problem in a multidimensional diffusion model, using the same HJBI approach as \citet{Oksendal2008HJBI} and \citet{trybula2019continuous}. {Recently, there are also works studying MMV investment problems when the trading strategy is constrained, see e.g., \citet{shen2022cone}, \citet{hu2023constrained} and \citet{du2023monotone}. They all conclude that the solutions to MMV and MV problems coincide when the asset is continuous.}

As the method and result of \citet{SL2020note} can only be applied to continuous markets, a crucial question is whether the conclusion still holds true in any L\'{e}vy market. It is also an open question raised by \citet{trybula2019continuous} for further study. Thus, to find the relationship between the two optimal strategies under MMV and MV preferences in the L\'{e}vy market, we consider a jump diffusion financial market. Besides, an untradable stochastic factor is also included in our model for two reasons. First, if we follow the model setup of \citet{trybula2019continuous}, then we can compare our result with the former work and the difference is the pure effect of the jump risk. Second, markets with stochastic factors are incomplete. \citet{BG2015cashflow} propose a detailed definition of the free cash flow stream and prove that there is no free cash flow stream in any complete and arbitrage-free market. That is why the optimal strategy of an MV investor would never escape the monotone domain of MV preferences in such market. Therefore, it is necessary to consider an incomplete market in our research. It is important to note that the jump diffusion part of our model can also lead to the incomplete property of the market, depending on the specific coefficient setting of the jump strength.

From another perspective, we can regard the MMV preference as a robust utility, which is considered in many robust optimization problems, e.g., \citet{GS1989maxmin}, \citet{zawisza2010robust} and  \citet{Oksendal2014modeluncertainty}. In a jump diffusion model, \citet{Oksendal2008HJBI} consider such kind of max-min problem as a stochastic differential game. An interesting characterization on the set of ambiguity probability measures is proved in \citet{LG2021RAIRO-OR}, which helps us understand the range of ambiguity probability measures of MMV preferences even in the jump diffusion market. We will discuss more about this in Remark \ref{MMV Q define}. In this paper, we mainly refer to \citet{Oksendal2008HJBI} and take a similar dynamic programming method to solve the MMV optimal investment problem in our market model.

Beside the non-monotonicity, another huge drawback of MV preferences is the time-inconsistency. The main reason for the time-inconsistency is the variance term in MV preferences, which yields the in-applicability of the dynamic programming method. \citet{BTMA2010generaltimeincontheory} propose a general theory of Markovian time inconsistent stochastic control problems, in which three approaches are raised to fix the time-inconsistency and find the optimal control. The first is to choose whatever seems to be optimal based on the investor's current preference, without knowing the fact that his preference may change over time. The second is to find a pre-commitment solution, which means that the investor only considers the optimal problem once at the initial time and keeps taking that strategy all the time. All of the literature on MV preferences mentioned above actually uses the pre-commitment approach. The third is to define a new equilibrium solution and use game theory terms to formulate the problem. For example, \citet{BTMA2010generaltimeincontheory} find an equilibrium solution to the general utility form that includes MV preferences. \citet{ZLL2013IMEmvtimeconsistent} derive the optimal strategy of reinsurance and investment that holds time-consistency for MV insurers. All three approaches can be reasonable in different circumstances, depending on the goal of the study. In our paper, similar to \citet{trybula2019continuous}, we use the pre-commitment approach, which is the best situation for us to analyze the relationship between MMV and MV problems. As we find that the optimal strategy under MMV preferences depend on the investor's wealth and the market state, MMV preferences also have the time-inconsistency property as MV preferences. 

The main contributions of our work are as follows. We refine the optimal investment problems with respect to MMV and MV preferences in a jump diffusion and stochastic factor model. To the best of our knowledge, this is the first time that these two preferences are considered in such a financial market model with jumps. {We obtain the explicit optimal strategies of MMV and MV problems under certain assumptions (Conditions of Proposition \ref{pde condition} and Assumption \ref{assumption jumpsize}). We prove that the two optimal strategies and value functions coincide if and only if Assumption \ref{assumption jumpsize} holds. Thus, the discontinuity of the market model may cause investors with MMV preferences to choose a different strategy from those with MV preferences. We answer the open question proposed by \citet{trybula2019continuous} and extend the conclusion of \citet{SL2020note}. This phenomenon is interesting because it contradicts the consistency in continuous markets, where the two preferences are shown to yield the same result. Besides, under Assumption \ref{assumption jumpsize}, we prove that the achieved maximum utilities of MMV and MV preferences are also the same, which is somehow ignored by \citet{trybula2019continuous}. Theoretically, Assumption \ref{assumption jumpsize} may not hold, leading to the inconsistency of MMV and MV problems. Nevertheless, we also provide some empirical evidences to illustrate the validity of Assumption \ref{assumption jumpsize} in real financial market.}

More specifically, {we introduce the Dol\'eans-Dade exponentials to replace the target probability measures of MMV preferences.} In the MMV problem, we transform the problem into a stochastic differential game and use the dynamic programming method to obtain the HJBI equation. In the MV problem, we deal with the nonlinearity of the variance term first, and then obtain the HJB equation by combining the Lagrange multiplier and dynamic programming methods. To solve the HJBI and HJB equations, we obtain a new type of highly nonlinear parabolic partial differential equation. The equation cannot be simply transformed to a linear parabolic equation by the traditional Hopf-Cole transformation method used in \citet{zariphopoulou2001unhedgeable}, \citet{zawisza2010robust} and \citet{trybula2019continuous}. Instead, we refer to a new type existence and uniqueness theorem on parabolic equations from \citet{parabolicLadyzhenskaya} to prove the existence and uniqueness of the solutions to HJBI and HJB equations. {Then, we propose a jump diffusion type Feynman-Kac formula to prove the consistency of the two optimal strategies of MMV and MV problems under Assumption \ref{assumption jumpsize}. Finally, we prove that the optimal strategies and value functions of MMV and MV problem are the same if and only if Assumption \ref{assumption jumpsize} holds. We also derive the efficient frontier of the investment and analyze the economic impact caused by the jump diffusion part of the risky asset.}

The remainder of this paper is organized as follows. Section \ref{section model} introduces the market model and monotone mean-variance preferences. Section \ref{section MMV} solves the optimal investment problem with respect to MMV preferences. We formulate the HJBI equation and find the corresponding solution. Section \ref{section mv} solves the optimal investment problem with respect to MV preferences. {In Section \ref{section relationship}, we compare the two optimal strategies and value functions and provide the sufficient and necessary condition for them to be consistent. Besides, we derive the efficient frontier and make some economic explanations. We also provide empirical evidences to demonstrate the validity of the assumption in real financial market.} Finally, Section \ref{section conclusion} concludes the paper. 

\section{Market model and MMV preference.}\label{section model}

We consider a filtered probability space $\(\Omega;\F;\{\F_t\}_{t\geq 0};\p\)$ satisfying the usual conditions. $W^1=\{W^1_t,t\in[0,T]\}$ and $W^2=\{W^2_t,t\in[0,T]\}$ are Brownian motions. {${N}$ is a Poisson random measure on $\([0,T]\times \R,\mathcal{B}([0,T]\times \R)\)$ with respect to a L\'evy process $L$. Filtration $\{\F_t\}_{t\geq 0}$ is generated by $W^1$, $W^2$ and $N$.} We assume that $W^1$, $W^2$ and $N$ are independent. Further more, we denote $\ud t\nu(\ud p)$ as the compensator (or intensity measure), and $\tilde{N}$ as the corresponding compensated martingale measure of $N$, i.e., $\tilde{N}\(\ud t,\ud p\)=N\(\ud t,\ud p\)-\ud t\nu(\ud p)$. We assume that $\nu(\R)< +\infty$. For more information about the Poisson random measure, see \citet{Protter2005} and \citet{tankov2003jumpmodeling}. 

Suppose that there are two investment assets. The price process of the risk-free asset $B=\{B_t,t\in[0,T]\}$ is given by 
\begin{equation}
	\ud B_t=B_t r\ud t,\ B_0=s_0.
\end{equation}
And the risky asset $S=\{S_t,t\in[0,T]\}$ follows a L\'{e}vy-It\"{o} process, satisfying the following stochastic differential equation (SDE):
\begin{equation}\label{dynamic of S}
	\ud S_t=S_{t-}\[\mu\(Z_t\)\ud t+\sigma\(Z_t\)\ud W^1_t+\int_{\mathbb{R}}\gamma\(Z_{t},p\) \tilde{N}\(\ud t,\ud p\)\],\ S_0=s,
\end{equation}
where the stochastic factor $Z=\{Z_t,t\in[0,T]\}$ follows SDE:
\begin{equation}\label{dynamic of Z}
	\ud Z_t=a\(Z_t\)\ud t+b\(Z_t\)\(\rho_W\ud W^1_t+\bar{\rho}_W\ud W^2_t\),\ Z_0=z. 
\end{equation}
Besides, we assume $\gamma> -1$ with $\int_{\R}\gamma(z,p)^2\nu(\ud p)<\infty$, coefficients $\mu,\sigma,a$ and $b$ are continuous functions, $\gamma$ is continuous with respect to $z$, and they satisfy all the required regularity conditions to ensure that SDEs \eqref{dynamic of S} and \eqref{dynamic of Z} have unique strong solutions. Let $\rho_W$ represent the correlation of Brownian motion parts of $S$ and $Z$ and $\rho_W^2+\bar{\rho}_W^2=1$.

In this paper, we compare the optimal investment strategies of investors under MMV preferences and MV preferences. Now, we introduce the MMV preferences first. Based on \citet{maccheroni2009portfolio}, the class of MMV preferences is defined by
$$V_\theta (X)=\inf_{\Q\in \mathcal{Q}_0}\{\E^{\Q}[X]+\frac{1}{2\theta}C\(\Q | \p\)\},$$
$$C\(\Q | \p\)=\E^{\p}\[\(\frac{\ud \Q}{\ud \p}\)^2\]-1,$$
where $$\mathcal{Q}_0=\left\{\Q:\Q\ll \p,\ \E^{\p}\[\(\frac{\ud \Q}{\ud \p}\)^2\]<+\infty \right\}.$$

Moreover, MMV preferences can be rewritten as $$V_\theta(X)=-\Lambda_\theta(X)-\frac{1}{2\theta},$$
where $$\Lambda_\theta(X)=\sup_{\Q\in \mathcal{Q}_0}\E^{\Q}\[-X-\frac{1}{2\theta}\frac{\ud \Q}{\ud \p}\].$$ 
For simplicity, we denote $\Lambda_{\theta}$ as the objective function used in the MMV investment problem.

{ Naturally, we consider the set of all the Dol\'eans-Dade exponential measures
	$$\mathcal{Q}=\left\{\Q^{\eta,\phi} \ll \p : \frac{\ud \Q^{\eta,\phi}}{\ud \p}=\mathcal{E}\(\int \eta_{1,t}\ud W^1_t+\int \eta_{2,t}\ud W^2_t+\int \int_{\mathbb{R}}\phi_{t}(p) \tilde{N}\(\ud t,\ud p\)\)_T,(\eta,\phi)\in\mathcal{M} \right\},$$
	to replace $\mathcal{Q}_0$, where the range of $(\eta,\phi)=(\eta_{1,t},\eta_{2,t},\phi_t)_{0\leq t\leq T}$ is given by
	$$\mathcal{M}=\left\{(\eta,\phi) \text{ is predictable}:\E^{\p}\[\(\frac{\ud \Q^{\eta,\phi}}{\ud \p}\)^2\]<\infty,\ \E^{\p}\[\frac{\ud \Q^{\eta,\phi}}{\ud \p}\]=1,\frac{\ud \Q^{\eta,\phi}}{\ud \p}\geq 0\right\}.$$
	For $\Q\in\mathcal{Q}$, we just restrict it to satisfy $\Q\ll\p$ and $\frac{\ud \Q}{\ud \p}\geq 0$ rather than $\Q\sim\p$ as in \citet{trybula2019continuous}. It is a notable difference, because we cannot ensure the target $\Q$ to be equivalent to $\p$ when the market faces jump risks. For the aim of clarity, we state the definition of the Dol\'eans-Dade exponential and stress one of its important properties here. 
	\begin{definition}\label{def exponential}
		For a semimartingale $M$ with $M_0=0$, the stochastic exponential (Dol\'{e}ans-Dade exponential) of $M$, written as $\mathcal{E}(M)$, is the unique semimartingale $X$  satisfying $$X_t=1+\int_{0}^{t}X_{s-}\ud M_s.$$ 
		Detailedly, $X$ can be represented by 
		\begin{equation}\label{eq repre of exponential}
			X_t=\exp\left\{M_t-\frac{1}{2}\[M,M\]^c_t\right\}\prod_{0<s\leq t}(1+\Delta X_s)\exp\big\{-\Delta X_s\big\}, 
		\end{equation}
		where the infinite product converges and is c\`adl\`ag with finite variation.
	\end{definition}
	
	For more information about the Dol\'{e}ans-Dade exponential, see \citet{Protter2005}, Section 8 of Chapter \uppercase\expandafter{\romannumeral2}. Furthermore, similar to the well-known Kazamaki's and Novikov's criteria, we need the following criterion to determine whether a Dol\'{e}ans-Dade exponential is a martingale and then $\Q^{\eta,\phi}$ to be a true probability measure. 
	
	\begin{lemma}[M\'{e}min's criterion for exponential martingale]\label{lemma memin}
		Let $M$ be a local martingale and $M^0_t=\sum_{s\leq t}{\Delta M_s\one_{|\Delta M_s|\geq \frac{1}{2}}}$. Let $M^1_t$ be the compensator of $M^0_t$: $M^1_t=M^0_t-\tilde{M^0_t}$, and $M^2_t=M_t-M^1_t$.
		If the compensator of process $$\[M^c,M^c\]_t+\sum_{0<s\leq t}|\Delta M^1_t|+\sum_{0<s\leq t}\(\Delta M^2_t\)^2$$ is bounded, then $\mathcal{E}(M)$ is a uniformly integrable martingale.
	\end{lemma}
	\proof{Proof.}
	See Theorem \uppercase\expandafter{\romannumeral3}-3 in \citet{memin2006decompositions} or Exercise 13 of Chapter \uppercase\expandafter{\romannumeral5} in \citet{Protter2005}.
	\hfill
	\endproof }
\vskip 5pt
Now, we introduce a new process $Y^{\eta,\phi}=\{Y^{\eta,\phi}_t,t\in[0,T]\}$ as follows:
\begin{equation} \label{eq dynamic of Y}
	\ud Y_s^{\eta,\phi}=Y_{s-}^{\eta,\phi}\[\eta_{1,s}\ud W^1_s+\eta_{2,s}\ud W^2_s+ \int_{\mathbb{R}}\phi_{s}(p) \tilde{N}\(\ud s,\ud p\)\],\ Y_t=y.
\end{equation}
{Based on the uniqueness of SDE's solution, we have $Y_T^{\eta,\phi}=y\frac{\ud \Q^{\eta,\phi}}{\ud \p}$ when starting with $t = 0$ and $Y_0=y$.} If we set $y=\frac{1}{2\theta}$, the objective function restricted by the maximization on $\mathcal{Q}$ becomes $$\Lambda_\theta(X)=\sup_{\Q^{\eta,\phi}\in \mathcal{Q}}\E^{\Q^{\eta,\phi}}\[-X-Y_T^{\eta,\phi}\].$$
{ Notably, when $y\neq 0$, $\frac{Y_T^{\eta,\phi}}{y} = \frac{\ud \Q^{\eta,\phi}}{\ud \p}$ should be greater than 0. Thus, not all the predictable $(\eta,\phi)$ meet the requirement. We consider a series of admissible sets of $(\eta,\phi)$:
	$$\begin{aligned}
		\mathcal{M}^t=\Big\{(\eta,\phi)=&(\eta_{1,s},\eta_{2,s},\phi_s)_{t\leq s\leq T} \text{ is predictable}: \text{The solution of Equation \eqref{eq dynamic of Y} satisfies} \\ & \E^{\p}\[(Y_T^{\eta,\phi})^2\]\leq +\infty,\E^{\p}\[Y_T^{\eta,\phi}\]=1,Y_T^{\eta,\phi}\geq 0,\text{ for } \forall y>0\Big\}, \ 0 \leq t \leq T.
	\end{aligned}
	$$
	Denote $\mathcal{Q}^t$ as the set of all the Dol\'eans-Dade exponential measures with respect to $\mathcal{M}^t$. It is trivial to see that $\mathcal{M}^0=\mathcal{M}$ and $\mathcal{Q}^0=\mathcal{Q}$.} Further, we denote $\E^{\Q^{\eta,\phi}}$ by $\E^{\eta,\phi}$ for simplicity. In the following discussion, we restrict the maximization (minimization) over $\Q$ on $\mathcal{Q}^t$ in the definition of $\Lambda_{\theta}$ ($V_{\theta}$), unless otherwise noted. See Remark \ref{MMV Q define} for more detailed explanations.
\begin{remark}\label{MMV Q define}
	In fact, the Dol\'eans-Dade exponential may not exhaust all the probability measure $\Q$ in $\mathcal{Q}_0$, which means $\mathcal{Q}^t \subset \mathcal{Q}_0$. Two reasons support us to consider $\mathcal{Q}^t$ as the target probability set. { First, we will prove in Section \ref{section relationship} that under certain conditions, the optimal control and value function of the MMV problem would be consistent no matter $\mathcal{Q}^t$ or $\mathcal{Q}_0$ is considered. Second, \citet{LG2021RAIRO-OR} propose that any $\Q$ in $\mathcal{Q}_0$ is indeed a Dol\'eans-Dade exponential if we consider $\{\F_t\}_{t\geq 0}$ as the filtration generated by Brownian motions, and $\Q$ as all the Dol\'eans-Dade exponential measures generated by Brownian motions. It inspires us to consider such form of probability measures even in the jump diffusion financial market.}
\end{remark}

\section{MMV problem.} \label{section MMV}

In this section, we solve the optimal investment problem under MMV preferences. We treat the max-min optimization problem as a stochastic differential game between an investor and a potential market player. We use the dynamic programming method to get the HJBI equation and propose a verification theorem. The optimal investment strategy is then obtained by solving the HJBI equation and verifying the verification theorem.

First, we introduce the investor's wealth process. Let $\bar{\pi}=\{\bar{\pi}_t,t\in[0,T]\}$ be the process of investment amount in the risky asset. Then, the wealth process $\bar{X}$ satisfies the following SDE:
\begin{equation*}
	\left\{\begin{aligned}
		&\ud \bar{X}_t^{\bar{\pi}}=\bar{X}_t^{\bar{\pi}}r\ud t+\bar{\pi}_{t}\(\mu\(Z_t\)-r\) \ud t+\bar{\pi}_t\sigma\(Z_t\)\ud W^1_t+\bar{\pi}_{t}\int_{\mathbb{R}}\gamma\(Z_{t},p\) \tilde{N}\(\ud t,\ud p\),\\ 
		&X_0=x.
	\end{aligned}\right.
\end{equation*}
Similar to \citet{trybula2019continuous}, we consider $\pi=\{e^{r(T-t)}\bar{\pi}_t,t\in[0,T]\}$ and $X^{{\pi}}=\{e^{r(T-t)}\bar{X}_t^{\bar{\pi}},t\in[0,T]\}$ as the forward value processes of the investment amount and wealth for simplicity. Then, the wealth process of the self-financed investor satisfies:
\begin{equation}\label{dynamic of Xt}
	\left\{\begin{aligned}
		&\ud X_t^{\pi}=\pi_{t}\(\mu\(Z_t\)-r\) \ud t+\pi_t\sigma\(Z_t\)\ud W^1_t+\pi_{t}\int_{\mathbb{R}}\gamma\(Z_{t},p\) \tilde{N}\(\ud t,\ud p\),\\ 
		&X_0=x.
	\end{aligned}\right.
\end{equation}

{\begin{definition}
		An investment strategy $\pi=\{\pi_s\}_{t\leq s\leq T}$ is admissible $(\pi\in \mathcal{A}_{x,y,z,t})$ if and only if
		\begin{enumerate}[(\romannumeral1)]
			\item $\pi$ is predictable.
			\item $\pi$ is such that SDE (\ref{dynamic of Xt}) has a unique solution satisfying $$\E^{\eta,\phi}_{x,y,z,t}\[\sup_{t\leq s\leq T}|X_s^{\pi}|^2\]<\infty,\ \text{for} \ \forall(\eta,\phi)\in \mathcal{M}^t.$$
		\end{enumerate}
\end{definition}}

Next, we define the objective function as 
$$J^{\pi,\eta,\phi}(x,y,z,t)=\E^{\eta,\phi}_{x,y,z,t}\[-X_T^{\pi}-Y_T^{\eta,\phi}\].$$ 

\begin{problem}\label{MMV problem}
	Solve the stochastic control problem: 
	\begin{equation}
		V(x,y,z,t)=\min_{\pi} \sup_{\eta,\phi}J^{\pi,\eta,\phi}(x,y,z,t).
	\end{equation}
	If $\(\pi^*,\eta^*,\phi^*\) \in \mathcal{A}_{x,y,z,t}\times \mathcal{M}^t$ satisfies $V(x,y,z,t)=J^{\pi^*,\eta^*,\phi^*}(x,y,z,t)$, then we call $\(\pi^*,\eta^*,\phi^*\)$ the optimal control of Problem \ref{MMV problem}, $V(x,y,z,t)$ the value function of Problem \ref{MMV problem} and $\Q^{\eta^*,\phi^*}$ the optimal measure of MMV preferences. 
\end{problem}
\vskip 5pt

\begin{remark}\label{MMV Girsanov}
	Setting $y=\frac{1}{2\theta}$, we can get the solution to the MMV investment problem. In the stochastic control problem, we allow it to be a changeable value for the application of the dynamic programming method. { Besides, as we need to take the expectation under a new probability measure in the objective function, we use a general form of the Girsanov Theorem which does not require $\Q^{\eta,\phi}\sim \p$ but only $\Q^{\eta,\phi}\ll \p$ (see, \citet{chan1999pricing} and Section 3 of Chapter 3 in \citet{jacod2013limit}), to rewrite the dynamics of $X,Y$ and $Z$ under the probability measure $\Q^{\eta,\phi}$.}
\end{remark}

Based on the Girsanov Theorem, $W^{\eta_1}_t=W^1_t-\int_{0}^{t}\eta_{1,s}\ud s$ and $W^{\eta_2}_t=W^2_t-\int_{0}^{t}\eta_{2,s}\ud s$ are Brownian motions under $\Q^{\eta,\phi}$. For the jump part, the compensator of $N\(\ud t,\ud p\)$ becomes $(1+\phi_{t}(p))\nu(\ud p)\ud t$, i.e., $\tilde{N}^{\phi}\(\ud t,\ud p\)=\tilde{N}\(\ud t,\ud p\)-\phi_{t}(p)\nu(\ud p)dt$ is compensated Poisson random measure under $\Q^{\eta,\phi}$. Then, the dynamics under $\Q^{\eta,\phi}$ become:
\begin{equation}
	\left\{\begin{aligned}
		\ud X_t^{\pi}=&\pi_{t}\[\mu\(Z_t\)-r+\sigma\(Z_t\)\eta_{1,t}-\int_{\mathbb{R}} \gamma\(Z_{t},p\)\nu(\ud p)  \] \ud t\\
		&+\pi_t\sigma\(Z_t\)\ud W^{\eta_1}_t+\pi_{t}\int_{\mathbb{R}}{\gamma\(Z_{t},p\)}N\(\ud t,\ud p\),\\ 
		\ud Y_t^{\eta,\phi}=&Y_t^{\eta,\phi}\[\eta_{1,t}^2+\eta_{2,t}^2-\int_{\mathbb{R}} \phi_{t}(p)\nu(\ud p)\]\ud t\\
		&+Y_{t-}^{\eta,\phi}\[\eta_{1,t}\ud W^{\eta_1}_t+\eta_{2,t}\ud W^{\eta_2}_t+ \int_{\mathbb{R}}\phi_{t}(p) N\(\ud t,\ud p\)\],\\ 
		\ud Z_t=&\big[a\(Z_t\)+b\(Z_t\)\(\rho_W\eta_{1,t}+\bar{\rho}_W\eta_{2,t}\)\big]\ud t +b\(Z_t\)\(\rho_W\ud W^1_t+\bar{\rho}_W\ud W^2_t\).
	\end{aligned}\right.
\end{equation}

Referring to \citet{Oksendal2007ch3}, the generator $\mathscr{A}$ of $(X,Y,Z)$ is:
\begin{equation}\label{MMV generator}
	\begin{aligned}
		\mathscr{A}^{\pi,\eta,\phi}V(x,y,z,t)
		=&V_t+V_x\pi\[\mu\(z\)-r+\sigma\(z\)\eta_{1}-\int_{\mathbb{R}} \gamma\(z,p\)\nu(\ud p) \]\\
		&+V_yy\[\eta_{1}^2+\eta_{2}^2-\int_{\mathbb{R}} \phi(p)\nu(\ud p)\]+V_z\[a\(z\)+b\(z\)\(\rho_W\eta_{1}+\bar{\rho}_W\eta_{2}\)\]\\
		&+\frac{1}{2}V_{xx}\pi^2\sigma(z)^2+\frac{1}{2}V_{yy}y^2\(\eta_{1}^2+\eta_{2}^2\)+\frac{1}{2}V_{zz}b(z)^2\\
		&+V_{xy}\pi\sigma(z)y\eta_1+V_{xz}\pi\sigma(z)\rho_Wb(z)+V_{yz}yb(z)\(\rho_W\eta_{1}+\bar{\rho}_W\eta_{2}\)\\
		&+\int_{\mathbb{R}}\[V(x+\pi\gamma(z,p),y+y\phi(p),z,t)-V(x,y,z,t)\]\(1+\phi(p)\)\nu(\ud p).
	\end{aligned}
\end{equation}
See \citet{Oksendal2007ch1} for more information about the Girsanov theorem and the generator. 

\subsection{Verification theorem.}
To solve Problem \ref{MMV problem}, we propose and prove the following verification theorem in this subsection. For simplicity, we slightly abuse the notation $\pi,\eta$ or $\phi$ to represent a process or a real number in $\R$.

\begin{theorem}\label{MMV verification}
	Suppose that there exist a function $$V\in C^{2,2,2,1}\(\R\times\R^+\times\R\times[0,T)\)\cap C\(\R\times\R^+\times\R\times[0,T]\)$$ and a control $$(\pi^*,\eta^*,\phi^*)\in \mathcal{A}_{x,y,z,t}\times \mathcal{M}^t$$ such that 
	\begin{enumerate}[(\romannumeral1)]
		\item $\mathscr{A}^{\pi^*,\eta,\phi}V(x,y,z,t)\leq 0,\ \forall \(x,y,z,t\) \in \R \times \R^+ \times \R \times [0,T),\ \forall (\eta,\phi) \in \R^2.$
		\item $\mathscr{A}^{\pi,\eta^*,\phi^*}V(x,y,z,t)\geq 0,\ \forall \(x,y,z,t\) \in \R \times \R^+ \times \R \times [0,T),\ \forall \pi \in \R.$
		\item $\mathscr{A}^{\pi^*,\eta^*,\phi^*}V(x,y,z,t)= 0,\ \forall\(x,y,z,t\) \in \R \times \R^+ \times \R \times [0,T).$
		\item $V(x,y,z,T)=-x-y\ a.s.,\  \forall\(x,y,z\) \in \R \times \R^+ \times \R .$
		\item $\E^{\eta,\phi}_{x,y,z,t}\[\sup_{t\leq s\leq T}|V(X_s^{\pi},Y_s^{\eta,\phi},Z_s,s))|\]<+\infty,\ \text{for} \ \forall\(x,y,z,t\) \in \R \times \R^+ \times \R \times [0,T), \forall (\pi,\eta,\phi)\in \mathcal{A}_{x,y,z,t}\times \mathcal{M}^t$.
	\end{enumerate}
	Then, we have $$J^{\pi^*,\eta,\phi}(x,y,z,t)\leq V(x,y,z,t)\leq J^{\pi,\eta^*,\phi^*}(x,y,z,t),\ \text{for} \ \forall \pi\in\mathcal{A}_{x,y,z,t},(\eta,\phi)\in\mathcal{M}^t,$$
	$$V(x,y,z,t)=J^{\pi^*,\eta^*,\phi^*}(x,y,z,t)=\min_{\pi} \sup_{\eta,\phi}J^{\pi,\eta,\phi}(x,y,z,t),$$
	$\(\pi^*,\eta^*,\phi^*\)$ is the optimal control and $V(x,y,z,t)$ is the value function of Problem \ref{MMV problem}. 
\end{theorem}

\proof{Proof.}
For $\forall \(\pi,\eta,\phi\)\in \mathcal{A}_{x,y,z,t}\times \mathcal{M}^t$, we define a stopping time $\tau_N=T\wedge \inf\{t\geq 0: |X_t^{\pi}|\geq N,|Y_t^{\eta,\phi}|\geq N,|Z_t|\geq N\}$. { Using  Ito's formula (see, e.g., \citet{Protter2005}), we have 
	\begin{equation}\label{dynkin formula}
		\E_{x,y,z,t}^{\eta,\phi}\[V(X_{\tau_N}^{\pi},Y_{\tau_N}^{\eta,\phi},Z_{\tau_N},\tau_N)\]=V(x,y,z,t)+\E_{x,y,z,t}^{\eta,\phi}\[\int_{t}^{\tau_N} \mathscr{A}^{\pi,\eta,\phi}V(X_{s-}^{\pi},Y_{s-}^{\eta,\phi},Z_s,s)\ud s\].
\end{equation}}

First, we take $\(\pi^*,\eta,\phi\)$ into Equation \eqref{dynkin formula}. Using Condition $(\romannumeral1)$, we have $$\E^{\eta,\phi}_{x,y,z,t}\[V(X_{\tau_N}^{\pi^*},Y_{\tau_N}^{\eta,\phi},Z_{\tau_N},\tau_N)\]\leq V(x,y,z,t).$$ 
Using Conditions $(\romannumeral4)$, $(\romannumeral5)$ and  Dominated convergence theorem,  letting $N\to \infty$, we obtain $$J^{\pi^*,\eta,\phi}(x,y,z,t)= \E^{\eta,\phi}_{x,y,z,t}\[-X_T^{\pi^*}-Y_T^{\eta,\phi}\]\leq V(x,y,z,t).$$  Thus, 
\begin{equation}\label{inequality 1}
	\min_{\pi} \sup_{\eta,\phi}J^{\pi,\eta,\phi}(x,y,z,t)\leq  \sup_{\eta,\phi}J^{\pi^*,\eta,\phi}(x,y,z,t)\leq V(x,y,z,t).
\end{equation}

Second, we replace $\(\pi,\eta^*,\phi^*\)$ into Equation \eqref{dynkin formula} and use the same way to obtain $$V(x,y,z,t)\leq  \E^{\eta^*,\phi^*}_{x,y,z,t}\[-X_T^{\pi}-Y_T^{\eta^*,\phi^*}\]=J^{\pi,\eta^*,\phi^*}(x,y,z,t).$$
Thus, 
\begin{equation}\label{inequality 2}
	V(x,y,z,t)\leq \min_{\pi} J^{\pi,\eta^*,\phi^*}(x,y,z,t)\leq \min_{\pi} \sup_{\eta,\phi}J^{\pi,\eta,\phi}(x,y,z,t).
\end{equation}
Combining Inequalities \eqref{inequality 1} and \eqref{inequality 2} yields $$V(x,y,z,t)=\min_{\pi} \sup_{\eta,\phi}J^{\pi,\eta,\phi}(x,y,z,t).$$

Finally, we take $\(\pi^*,\eta^*,\phi^*\)$ into Equation \eqref{dynkin formula} and use the same way to conclude  $$V(x,y,z,t)=J^{\pi^*,\eta^*,\phi^*}(x,y,z,t).$$ \hfill
\endproof 

{\begin{remark}\label{rmk MMV HJBI}
		Consider the following equation related to Conditions $(\romannumeral1)\sim(\romannumeral4)$ :
		\begin{equation}\label{MMV HJBI}
			\left\{\begin{aligned}
				&\min_{\pi} \sup_{\eta,\phi}\mathscr{A}^{\pi,\eta,\phi}V(x,y,z,t)=0,\\
				&V(x,y,z,T)=-x-y.
			\end{aligned}\right.
		\end{equation}
		We call Equation \eqref{MMV HJBI} the Hamilton–Jacobi–Bellman–Isaacs (HJBI) equation of Problem \ref{MMV problem} with the boundary condition. Next, we guess and verify the possible solution of this equation, and prove that it satisfies Conditions of Theorem \ref{MMV verification}. 
\end{remark}}   

\subsection{Solution to Problem \ref{MMV problem}.}

Based on the boundary condition $V(x,y,z,T)=-x-y$, and $X,Y$ are of the linear form in $J^{\pi,\eta,\phi}(x,y,z,t)$, we guess the value function has the following form: 
\begin{equation}\label{MMV value function form}
	V(x,y,z,t)=-x+G(z,t)y,
\end{equation}
where $G(z,T)=-1$. Substituting it into the generator \eqref{MMV generator}, we obtain
\begin{equation}\label{MMV generator 2}
	\begin{aligned}
		\mathscr{A}^{\pi,\eta,\phi}V(x,y,z,t)
		=&yG_t-\pi\[\mu\(z\)-r+\sigma\(z\)\eta_{1}-\int_{\mathbb{R}} \gamma\(z,p\)\nu(\ud p) \]\\
		&+yG\[\eta_{1}^2+\eta_{2}^2-\int_{\mathbb{R}} \phi(p)\nu(\ud p)\]+yG_z\[a\(z\)+b\(z\)\(\rho_W\eta_{1}+\bar{\rho}_W\eta_{2}\)\]\\
		&+\frac{1}{2}yG_{zz}b(z)^2+G_{z}yb(z)\(\rho_W\eta_{1}+\bar{\rho}_W\eta_{2}\)\\
		&+\int_{\mathbb{R}}\[-\pi\gamma(z,p)+Gy\(1+\phi(p)\)-Gy\]\(1+\phi(p)\)\nu(\ud p).
	\end{aligned}
\end{equation}

The generator above is a quadratic function of $\eta_{1}$, $\eta_{2}$ and $\phi$. To solve Equation \eqref{MMV HJBI}, we first fix $\pi$ and find the maximum point over $\eta_{1}$, $\eta_{2}$ and $\phi$. The first order condition is
$$\left\{\begin{aligned}
	&2yG(z,t)\eta_1-\pi\sigma(z)+2yG_z(z,t)b(z){\rho}_W=0,\\
	&2yG(z,t)\eta_2+2yG_z(z,t)b(z)\bar{\rho}_W=0,\\
	&2yG(z,t)\phi(p)-\pi\gamma(z,p)=0.
\end{aligned}\right. $$
If $G(z,t)<0$, then the Hessian matrix would be negatively definite and we obtain the local maximum point 
\begin{equation}\label{z etawithpi}
	\begin{aligned}
		&\eta_1^*(\pi)=\frac{\pi\sigma(z)}{2yG(z,t)}-\frac{G_z(z,t)b(z)\rho_W}{G(z,t)},\\
		&\eta_2^*(\pi)=-\frac{G_z(z,t)b(z)\bar{\rho}_W}{G(z,t)},\\
		&\phi^*(\pi,p)=\frac{\pi\gamma(z,p)}{2yG(z,t)}.
	\end{aligned}
\end{equation}
Substituting the candidates $\eta^*(\pi)$ and $\phi^*(\pi)$ into the generator \eqref{MMV generator 2}, we obtain
\begin{equation}\label{MMV generator 3}
	\begin{aligned}
		\mathscr{A}^{\pi,\eta^*(\pi),\phi^*(\pi)}V(x,y,z,t)
		=&yG_t-\pi\[\mu\(z\)-r-\int_{\mathbb{R}} \gamma\(z,p\)\nu(\ud p) \]+yG_za\(z\)+\frac{1}{2}yG_{zz}b(z)^2\\
		&+\int_{\mathbb{R}}\(-\pi\gamma(z,p)\)\nu(\ud p)-\frac{\(-\pi\sigma(z)+2yG_zb(z){\rho}_W\)^2}{4yG}\\
		&-\frac{\(2yG_zb(z)\bar{\rho}_W\)^2}{4yG}-\int_{\R}\frac{\(-\pi\gamma(z,p)\)^2}{4yG}\nu\(\ud p\).
	\end{aligned}
\end{equation}

Again, the above equation is a quadratic function of $\pi$. To find the minimum point over $\pi$, we consider the first order condition
$$-\frac{\sigma(z)^2+\int_{\R}{\gamma(z,p)^2}\nu\(\ud p\)}{2yG(z,t)}\pi-\(\mu\(z\)-r\)+\frac{\sigma(z)b(z){\rho}_WG_z(z,t)}{G(z,t)}=0.$$
If $G(z,t)<0$, then the local minimum point is given by 
\begin{equation}\label{candidate pi}
	\pi^*=-\frac{2yG(z,t)}{\sigma(z)^2+\int_{\R}{\gamma(z,p)^2}\nu\(\ud p\)}\[\mu\(z\)-r-\frac{\sigma(z)b(z){\rho}_WG_z(z,t)}{G(z,t)}\].
\end{equation}
To simplify the notation, we define
$$\Sigma(z)={\sigma(z)^2+\int_{\R}{\gamma(z,p)^2}\nu\(\ud p\)}.$$

Substituting the candidate $\pi^*$ given by \eqref{candidate pi} into \eqref{z etawithpi}, we get the candidate optimal control
\begin{equation}\label{MMV optimal control}
	\begin{aligned}
		&\pi^*=-\frac{2yG(z,t)}{\Sigma(z)}\[\mu\(z\)-r-\frac{\sigma(z)b(z){\rho}_WG_z(z,t)}{G(z,t)}\],\\
		&\eta_{1}^*=-\frac{\sigma(z)}{\Sigma(z)}\(\mu\(z\)-r\)-{\rho}_Wb(z)\frac{\int_{\R}{\gamma(z,p)^2\nu\(\ud p\)}}{\Sigma(z)}\frac{G_z(z,t)}{G(z,t)},\\
		&\eta_2^*=-\frac{G_z(z,t)b(z)\bar{\rho}_W}{G(z,t)},\\
		&\phi^*(p)=-\frac{\gamma(z,p)}{\Sigma(z)}\[\mu\(z\)-r-\frac{\sigma(z)b(z){\rho}_WG_z(z,t)}{G(z,t)}\].
	\end{aligned}
\end{equation}
Substituting them into \eqref{MMV generator 3}, we obtain
$$\begin{aligned}
	\mathscr{A}^{\pi^*,\eta^*,\phi^*}V(x,y,z,t)
	=&yG_t+yG_z\[a\(z\)-2\rho_W\frac{\mu\(z\)-r}{\Sigma(z)}\sigma(z)b(z)\]\\
	&+\frac{1}{2}yG_{zz}b(z)^2+yG\frac{\(\mu\(z\)-r\)^2}{\sigma(z)^2}-y\frac{\bar{\rho}_W^2\sigma(z)^2+\int_{\R}{\gamma(z,p)^2\nu\(\ud p\)}}{\Sigma(z)}b(z)^2\frac{G_z^2}{G}.
\end{aligned}$$
Eliminating $y$, the HJBI Equation \eqref{MMV HJBI} is transformed to the following PDE with respect to $G$:
\begin{equation}\label{MMV pde G}
	\left\{\begin{aligned}
		&G_t+\alpha(z)G_z+\frac{1}{2}b(z)^2G_{zz}-\beta(z)b(z)^2\frac{G_z^2}{G}+\lambda(z)G=0,\\
		&G(z,T)=-1,
	\end{aligned}\right.
\end{equation}
where 
$$\left\{\begin{aligned}
	&\alpha(z)\triangleq a\(z\)-2\rho_W\frac{\mu\(z\)-r}{\Sigma(z)}\sigma(z)b(z),\\
	&\beta(z)\triangleq\frac{\bar{\rho}_W^2\sigma(z)^2+\int_{\R}{\gamma(z,p)^2\nu\(\ud p\)}}{\Sigma(z)}>0,\\
	&\lambda(z)\triangleq\frac{\(\mu\(z\)-r\)^2}{\Sigma(z)}>0.
\end{aligned}\right.$$

If Equation \eqref{MMV pde G} has a proper solution $G$, we may derive out the optimal strategy and the correspinding value function based on Theorem \ref{MMV verification}. Thus, we temporarily focus on the property of Equation \eqref{MMV pde G} and provide the existence and regularity condition of the solution. By taking the Hopf-Cole transformation $G=-e^F$ (see \citet{zawisza2010robust} and \citet{liu2017hopf} for more information about such technique), Equation \eqref{MMV pde G} becomes
\begin{equation}\label{MMV pde}
	\left\{\begin{aligned}
		&F_t+\alpha(z)F_z+\frac{1}{2}b(z)^2F_{zz}+\[\frac{1}{2}-\beta(z)\]b(z)^2F_z^2+\lambda(z)=0,\\
		&F(z,T)=0.
	\end{aligned}\right.
\end{equation}

\begin{remark}\label{rmk PDE}
	To the best of our knowledge, it is the first time to solve such type of nonlinear parabolic equations in portfolio selection problems, where the nonlinear coefficient $\frac{1}{2}-\beta(z)$ is not a constant. For the constant case, Equation \eqref{MMV pde G} or \eqref{MMV pde} can be directly transformed to a classical linear parabolic equation, see \citet{zariphopoulou2001unhedgeable}, \citet{zawisza2010robust} and \citet{trybula2019continuous}. Different from them, here we need the depiction of quasi-linear parabolic equations in divergence form given by \citet{parabolicLadyzhenskaya}.
	
	{In the following discussion, we say a function $u\in C^{2+\beta,1+\beta/2}(\Omega)$ if and only if
		$$\begin{aligned}
			&\sup_{(x,t)\in \Omega}|u|+\sup_{(x,t)\in \Omega}|u_t|+\sup_{(x,t)\in \Omega}|u_x|+\sup_{(x,t)\in \Omega}|u_{xx}|\\&+\sup_{(x,t),(x',t')\in \Omega}\frac{|u_{xx}(x,t)-u_{xx}(x',t')|}{d((x,t),(x',t'))^\beta}+\sup_{(x,t),(x',t')\in \Omega}\frac{|u_{t}(x,t)-u_{t}(x',t')|}{d((x,t),(x',t'))^\beta}<\infty.
		\end{aligned}$$
		Similarly, a function $u\in C^{2,1}(\Omega)$ if and only if 
		$\sup\limits_{(x,t)\in \Omega}\big(|u|+|u_t|+|u_x|+|u_{xx}|\big)<\infty.$ Obviously, $u\in C^{2+\beta,1+\beta/2}(\Omega)$ indicates that $u$ and $u_x$ are bounded.}
\end{remark}  

\begin{theorem}[Cauchy problem for quasi-linear parabolic equations in divergence form]\label{parabolic equation}
	Consider the divergence form of parabolic equations:
	\begin{equation}\label{cauchy problem}
		\left\{\begin{aligned}
			&-u_t+\(f(x,t,u,u_x)\)_x+g(x,t,u,u_x)=0,\ (x,t)\in H_T=\R\times(0,T],\\
			&u(x,0)=\varphi(x),\ x\in \R.
		\end{aligned}
		\right.
	\end{equation}				
	Suppose that $f$ and $g$ satisfy the following conditions, where $z$ represents $u$ and $p$ represents $u_x$: 
	\begin{enumerate}[(\romannumeral1)]\label{pde thm condition}
		\item Denote $h(x,t,z,p)=g(x,t,z,p)+\frac{\partial f}{\partial z}p+\frac{\partial f}{\partial x}$. For $(x,t)\in {H_T}$ and any $(z,p)$, we have $$\frac{\partial f}{\partial p}\geq0,\ zh(x,t,z,0)\leq C_1z^2+C_2.$$
		\item For $(x,t)\in \bar{Q_T}$, $|z|\leq M$ and any $p$, where ${Q_T}=[-A,A]\times(0,T]$ and $A$ is arbitrary, we have that functions $f(x,t,z,p)$, $g(x,t,z,p)$ are continuous and $f$ is differentiable with respect to $x,z$ and $p$, and satisfy
		$$ 0<\lambda\leq \frac{\partial f(x,t,z,p)}{\partial p}\leq \Lambda, \ \(|f|+|\frac{\partial f}{\partial z}|\)\(1+|p|\)+|\frac{\partial f}{\partial x}|+|g|\leq C_3\(1+|p|^2\). $$ 
		\item For $(x,t)\in \bar{Q_T}$, $|z|\leq M$, $|p|\leq M$, we have that $f,g,\frac{\partial f}{\partial p},\frac{\partial f}{\partial z}$, $\frac{\partial f}{\partial x}$ are continuous and $\beta,\beta/2,\beta,\beta$-order H$\ddot{o}$lder continuous with respect to $x,t,z,p$.
		\item For any bounded interval $I\subset \R$, $\varphi(x)\in C^{2+\beta}(I)$ and $\max_{\R}|\varphi(x)|<\infty$. 
	\end{enumerate}
	Then, the Cauchy problem \eqref{cauchy problem} has a bounded solution $u$ and $u\in C^{2+\beta,1+\beta/2}(\bar{Q_T})$ for any $Q_T$. Furthermore, if constants $\lambda,\Lambda$ and $C_3$ do not depend on $Q_T$, then $u\in C^{2+\beta,1+\beta/2}({H_T})$. 
	
	Additionally, suppose that $f$ and $g$ satisfy the following condition:
	\begin{enumerate}[(\romannumeral5)]
		\item $\frac{\partial f}{\partial p}$ and $h$ are differentiable with respect to $z,p$, and satisfy
		$$ \max_{(x,t)\in H_T,|z|,|p|\leq M}|\frac{\partial^2 f}{\partial z\partial p},\frac{\partial^2 f}{\partial p^2},\frac{\partial h}{\partial p}|\leq C_4(M), \ \min_{(x,t)\in H_T,|z|,|p|\leq M} |\frac{\partial h}{\partial z}|\geq -C_5(M). $$
	\end{enumerate}
	Then, the solution to the Cauchy problem \eqref{cauchy problem} is unique.
\end{theorem}

{ In the conclusion of the above theorem, the solution $u\in C^{2+\beta,1+\beta/2}({H_T})$ indicates that $u$ and $u_x$ are bounded, which is an important fact to be used later.} Besides, we only state the one-dimensional case in the above theorem. The conclusion of the multi-dimensional case is nearly the same, see \citet{parabolicLadyzhenskaya}. Now we are ready to state the conclusion about Equation \eqref{MMV pde}.

\begin{proposition}\label{pde condition}
	Suppose that $b$ is continuously differentiable, $a$, $\sigma$ and {$ m(z):=\int\gamma(z,p)^2\nu(\ud p)$ are differentiable.} $b^2$ is bounded below from 0 and $a,b,b'$ and $\lambda$ are bounded. Then, Equation \eqref{MMV pde} has a unique solution $F\in C^{2,1}(H_T)$. Furthermore, if $b'$ is $\beta$-H$\ddot{o}$lder continuous, then $F\in C^{2+\beta,1+\beta/2}(H_T)$.
\end{proposition}

\proof{Proof.}
To apply Theorem \ref{parabolic equation}, we need to transform Equation \eqref{MMV pde} to the divergence form. First, we consider the inverse time form of Equation \eqref{MMV pde}:
\begin{equation}\label{MMV pde3}
	\left\{\begin{aligned}
		&-F_t+\frac{1}{2}b(z)^2F_{zz}+\alpha(z)F_z+\(\frac{1}{2}-\beta(z)\)b(z)^2F_z^2+\lambda(z)=0,\\
		&F(z,0)=0.
	\end{aligned}\right.
\end{equation}
Then, by simple calculation and comparison of Equations \eqref{cauchy problem} and \eqref{MMV pde3}, we have
$$\begin{aligned}
	&f(x,t,z,p)=\frac{1}{2}b(x)^2p,\\
	&g(x,t,z,p)=\(\alpha(x)-b(x)b'(x)\)p+\frac{1}{2}b(x)^2\beta(x)p^2+\lambda(x),\\
	&h(x,t,z,p)=g(x,t,z,p)+b(x)b'(x)p=\alpha(x)p+\frac{1}{2}b(x)^2\beta(x)p^2+\lambda(x).
\end{aligned}
$$
Now, we show that $f,g$ and $h$ satisfy Conditions of Theorem \ref{parabolic equation}. As $$\frac{\partial f}{\partial p}=\frac{1}{2}b(x)^2\geq0,\ zh(x,t,z,0)=z\lambda(x)\leq C_1z^2+C_2,$$
Condition $(\romannumeral1)$ of Theorem \ref{parabolic equation} holds. Because
$$\begin{aligned}
	|\alpha(x)|\leq& |a\(x\)|+2|\rho_W\frac{\mu\(x\)-r}{\Sigma(x)}\sigma(x)b(x)|\\
	\leq&|a\(x\)|+|\rho_Wb(x)\[\lambda(x)+\frac{\sigma(x)^2}{\Sigma(x)}\]| \leq |a\(x\)|+|\rho_Wb(x)\[\lambda(x)+1\]|
\end{aligned}$$is bounded, and $0\leq \beta(x)\leq1,$
we have 
$$0<\lambda\leq \frac{1}{2}b(x)^2\leq \Lambda,$$
$$\begin{aligned}
	\(|f|+|\frac{\partial f}{\partial z}|\)\(1+|p|\)+|\frac{\partial f}{\partial x}|+|g| &=\frac{1}{2}b^2|p|(1+|p|)+|bb'||p|+|\(\alpha-bb'\)p+\frac{1}{2}b^2\beta p^2+\lambda|\\
	&\leq C_3\(1+|p|^2\),
\end{aligned}$$and the constant $C_3$ only depends on the bound of $a,b,b',\lambda$. While the continuity and differentiability of functions are easy to check,as such Condition $(\romannumeral2)$ of  Theorem \ref{parabolic equation} holds. 

To show Condition $(\romannumeral3)$, we have $\frac{\partial f}{\partial p}=\frac{1}{2}b^2,\frac{\partial f}{\partial z}=0,\frac{\partial f}{\partial x}=bb'p$. 
When $b'$ is continuous, the largest exponent of the H$\ddot{o}$lder continuity of $f,g,\frac{\partial f}{\partial p},\frac{\partial f}{\partial z},\frac{\partial f}{\partial x}$ is 0. When $b'$ is $\beta$-H$\ddot{o}$lder continuous, the largest exponent of the H$\ddot{o}$lder continuity is $\beta$. Thus, Condition $(\romannumeral3)$ holds. 

Condition $(\romannumeral4)$ is automatically satisfied because we have the boundary condition $\varphi(x)\equiv0$. As constants $C_1,C_2$ and $C_3$ in inequalities only depend on the bound of $a,b,b'$ and $\lambda$, by using Theorem \ref{parabolic equation}, Equation \eqref{MMV pde} has a solution $F\in C^{2,1}(H_T)$(or $C^{2+\beta,1+\beta/2}({H_T})$).

Furthermore, for Condition $(\romannumeral5)$, we have 
$$\max_{(x,t)\in H_T,|z|,|p|\leq M}|\frac{\partial^2 f}{\partial z\partial p},\frac{\partial^2 f}{\partial p^2},\frac{\partial h}{\partial p}|=\max_{(x,t)\in H_T,|z|,|p|\leq M}|\alpha(x)+\frac{1}{2}b^2(x)\beta(x)|\leq C_4,$$
$$\min_{(x,t)\in H_T,|z|,|p|\leq M} |\frac{\partial h}{\partial z}|=0.$$ Thus, we conclude that Equation \eqref{MMV pde} has a unique solution $F\in C^{2,1}(H_T)$ (or $C^{2+\beta,1+\beta/2}({H_T})$).		

\hfill 
\endproof 

Proposition \ref{pde condition} provides the existence of the solution to Equation \eqref{MMV pde} under certain assumptions. Next, we foucus on finding the optimal strategy and value function of Problem \ref{MMV problem}. Based on \eqref{MMV optimal control}, the candidate optimal control $\pi^*$ is relevant to the process $Y^{\eta^*,\phi^*}$. However, as it is untradable and cannot be instantaneously observed by the investor, such form of the optimal strategy is unsuitable in real investment environment. We provide the following relationship between $X^{\pi^*}$ and $Y^{\eta^*,\phi^*}$ under the candidate optimal control, which ensures that $\pi^*_t$ can be also represented by $X_t^{\pi^*}$.

\begin{theorem}\label{Y representation}
	Under the candidate optimal control \eqref{MMV optimal control}, $Y^{\eta^*,\phi^*}$ has the following representation:
	\begin{equation}\label{Y representation equation}
		2Y^{\eta^*,\phi^*}_sG(Z_s,s)=X^{\pi^*}_s-x+2yG(z,t),\ \forall s\in [t,T].
	\end{equation}
\end{theorem}
\proof{Proof.}

As both sides of Equation \eqref{Y representation equation} have the same initial value $2yG(z,t)$, we just need to check that they follow the same dynamic. Substituting the candidate optimal control \eqref{MMV optimal control} into Equations \eqref{dynamic of Z}$\sim$\eqref{dynamic of Xt}, we have
$$\begin{aligned}
	\ud X_t^{\pi^*}=&-\frac{2Y^{\eta^*,\phi^*}_{t-}G}{\Sigma}\[\mu-r-\frac{\sigma b{\rho}_WG_z}{G}\]\[\mu-r-\int_{\R}{\gamma}\nu\(\ud p\)\]\ud t\\
	&-\frac{2Y^{\eta^*,\phi^*}_{t-}G}{\Sigma}\[\mu-r-\frac{\sigma b{\rho}_WG_z}{G}\]\[\sigma\ud W_t^1+\int_{\R}{\gamma}N\(\ud t,\ud p\)\],\\
	&\\
	\ud Y^{\eta^*,\phi^*}_t=&Y^{\eta^*,\phi^*}_{t-}\[\frac{\int_{\R}\gamma\nu\(\ud p\)}{\Sigma}\(\mu-r\)-\frac{ b{\rho}_WG_z}{G}\frac{\sigma\int_{\R}\gamma\nu\(\ud p\)}{\Sigma}\]\ud t\\[4pt]
	&+Y^{\eta^*,\phi^*}_{t-}\[\(-\frac{\sigma}{\Sigma}\(\mu-r\)-{\rho}_Wb\frac{\int_{\R}{\gamma^2\nu\(\ud p\)}}{\Sigma}\frac{G_z}{G}\)\ud W_t^1-\frac{G_zb\bar{\rho}_W}{G}\ud W_t^2\]\\[4pt]
	&+Y^{\eta^*,\phi^*}_{t-}\[\int_{\R}-\frac{\gamma}{\Sigma}\(\mu-r\)+\frac{ b{\rho}_WG_z}{G}\frac{\sigma\gamma}{\Sigma}N\(\ud t,\ud p\)\],\\
	&\\
	\ud G(Z_t,t)=&\[\frac{2\rho_W\sigma b}{\Sigma}(\mu-r)G_z+b^2\frac{\bar{\rho}_W^2\sigma^2+\int_{\R}{\gamma^2\nu\(\ud p\)}}{\Sigma}\frac{G_z^2}{G}-\lambda G\]\ud t\\[4pt]
	&+G_zb\(\rho_W\ud W^1_t+\bar{\rho}_W\ud W^2_t\),
\end{aligned}$$
where we have omitted $Z_t,t,p$ in the notations of $G(Z_t,t),G_z(Z_t,t),b(Z_t),\sigma(Z_t),\gamma(Z_t,p)$ and $\Sigma(Z_t)$ for simplicity. Then, using  Ito's formula, we obtain
$$\begin{aligned}
	\ud Y^{\eta^*,\phi^*}_{t}G(Z_t,t)
	=&Y^{\eta^*,\phi^*}_{t-}\ud G(Z_t,t)+G(Z_t,t)\ud Y^{\eta^*,\phi^*}_t+\ud \[Y^{\eta^*,\phi^*}, G(Z,\cdot)\]_t\\[4pt]
	=&Y^{\eta^*,\phi^*}_{t-}\bigg\{G_z\frac{2\rho_W\sigma b}{\Sigma}(\mu-r)+\frac{G_z^2}{G}b^2\frac{\bar{\rho}_W^2\sigma^2+\int_{\R}{\gamma^2\nu\(\ud p\)}}{\Sigma}\\[4pt]
	&- G\lambda+G\[\frac{\int_{\R}\gamma\nu\(\ud p\)}{\Sigma}\(\mu-r\)-\frac{ b{\rho}_WG_z}{G}\frac{\sigma\int_{\R}\gamma\nu\(\ud p\)}{\Sigma}\]\\[4pt]
	&+G_zb{\rho}_W\[-\frac{\sigma}{\Sigma}\(\mu-r\)-{\rho}_Wb\frac{G_z}{G}\frac{\int_{\R}{\gamma^2\nu\(\ud p\)}}{\Sigma}\]+G_zb\bar{\rho}_W\(-\frac{G_zb\bar{\rho}_W}{G}\)\bigg\}\ud t\\[4pt]
	&+Y^{\eta^*,\phi^*}_{t-}\[G_zb\rho_W+G\(-\frac{\sigma}{\Sigma}\(\mu-r\)-{\rho}_Wb\frac{G_z}{G}\frac{\int_{\R}{\gamma^2\nu\(\ud p\)}}{\Sigma}\)\]\ud W_t^1\\[4pt]
	&+Y^{\eta^*,\phi^*}_{t-}\[G_zb\bar{\rho}_W+G\(-\frac{G_zb\bar{\rho}_W}{G}\)\]\ud W_t^2\\[4pt]
	&+GY^{\eta^*,\phi^*}_{t-}\int_{\R}\[-\frac{\gamma}{\Sigma}\(\mu-r\)+\frac{ b{\rho}_WG_z}{G}\frac{\sigma\gamma}{\Sigma}\]N\(\ud t,\ud p\)\\[4pt]
	=&-\frac{Y^{\eta^*,\phi^*}_{t-}G}{\Sigma}\(\mu-r-\frac{\sigma b{\rho}_WG_z}{G}\)\[\mu-r-\int_{\R}{\gamma}\nu\(\ud p\)\]\ud t\\[4pt]
	&-\frac{Y^{\eta^*,\phi^*}_{t-}G}{\Sigma}\(\mu-r-\frac{\sigma b{\rho}_WG_z}{G}\)\[\sigma\ud W_t^1+\int_{\R}{\gamma}N\(\ud t,\ud p\)\],
\end{aligned}$$
which is exactly $\ud X_t^{\pi^*}/2$.\hfill
\endproof 

\begin{remark}
	Theorem \ref{Y representation} reveals the time-inconsistency of MMV preferences. It is widely known that as MV preferences violate the law of iterated expectation and have the property of time-inconsistency when we solve the related control problem. Based on the \eqref{MMV optimal control} and \eqref{Y representation equation}, the representation of $\pi^*$ includes the initial value $x,y$ and $z$, which means that the candidate optimal investment strategy depends on the investor's wealth and the market state. Thus, MMV preferences also have the time-inconsistency property as the classical MV preferences. 
\end{remark}

{ Last, we prove that the candidate optimal control is indeed the optimal control. Before that, we need the following lemma to determine when a Dol\'eans-Dade exponential is a martingale. 
	
	\begin{lemma}\label{lem stochastic exponencial}
		Suppose that $\eta_1(z)$, $\eta_2(z)$ and $\phi(z,p)$ are bounded measurable functions. Further suppose  
		$$\phi(Z_{s\wedge \tau},\Delta L_{s\wedge \tau})\one_{\Delta L_{s\wedge \tau}\neq 0}\geq -1,\ \forall s\in(t,T], \ a.s.$$ where $\tau=\inf\{s> t:\phi(Z_s,\Delta L_s)\one_{\Delta L_{s}\neq 0}= -1\}$. Then, $(\eta,\phi)=(\eta_1(Z_s),\eta_2(Z_s),\phi(Z_s,p))_{t\leq s\leq T}$ is admissible and belongs to $\mathcal{M}^t$.
	\end{lemma}
	
	\proof{Proof.}
	Let $M_s=\int_{t}^{s}\eta_1(Z_u)\ud W^1_u+\int_{t}^{s}\eta_2(Z_u)\ud W^2_u+ \int_{t}^{s}\int_{\mathbb{R}}\phi(Z_u,p) \tilde{N}\(\ud u,\ud p\)$. It suffices to prove that $\e(M)$ is a martingale and $\e(M)\geq 0$. Because both the jump size and coefficients of Brownian motions of $M$ are bounded, it satisfies the Memin's criteria (Lemma \ref{lemma memin}) and implies that $\e(M)$ is a true martingale. To find that $\e(M)\geq 0$, we recall that $\e(M)$ has the representation \eqref{eq repre of exponential}. Then $\e(M)_T\geq 0$ is equivalent to $$\prod_{t<s\leq T}(1+\phi(Z_{s},\Delta L_{s})\one_{\Delta L_{s}\neq 0})\exp\{-\phi(Z_{s},\Delta L_{s})\one_{\Delta L_{s}\neq 0}\} \geq 0,$$ which is guaranteed by the assumption.
	\hfill
	\endproof 
	
	Based on Lemma \ref{lem stochastic exponencial}, we introduce a crucial assumption on the parameter setting of our financial market, which is closely related to the verification theorem of MMV problem and the equivalency of MMV and MV problems. It plays an important role in our further discussions.
	
	\begin{assumption}\label{assumption jumpsize}
		Define  stopping time $\tau=\inf\{s> t:\zeta(Z_s,s)\gamma(Z_s,\Delta L_s)\one_{\Delta L_{s}\neq 0}= -1\}$, where $L$ is the L\'evy process corresponding to the Poisson random measure $N$. We assume 
		$$\zeta(Z_{s\wedge \tau},s\wedge \tau)\gamma(Z_{s\wedge \tau},\Delta L_{s\wedge \tau})\one_{\Delta L_{s\wedge \tau}\neq 0}\geq -1,\ \forall s\in(t,T], \ a.s..$$
\end{assumption}}

Now, we are ready to state and prove the solution to the MMV Problem \ref{MMV problem}.
\begin{theorem}\label{MMV optimal control thm}
	Under Conditions of Proposition \ref{pde condition} and Assumption \ref{assumption jumpsize}, the candidate optimal control
	\begin{equation}\label{MMV optimal control final}
		\begin{aligned}
			\pi^*_s=&-\frac{2Y^{\eta^*,\phi^*}_{s-}G(Z_s,s)}{\Sigma(Z_s)}\[\mu\(Z_s\)-r-\frac{\sigma(Z_s)b(Z_s){\rho}_WG_z(Z_s,s)}{G(Z_s,s)}\],\\
			\eta_{1,s}^*=&-\frac{\sigma(Z_s)}{\Sigma(Z_s)}\(\mu\(Z_s\)-r\)-{\rho}_Wb(Z_s)\frac{\int_{\R}{\gamma(Z_s,p)^2\nu\(\ud p\)}}{\Sigma(Z_s)}\frac{G_z(Z_s,s)}{G(Z_s,s)},\\
			\eta_{2,s}^*=&-\frac{G_z(Z_s,s)b(Z_s)\bar{\rho}_W}{G(Z_s,s)},\\
			\phi_s^*(p)=&-\frac{\gamma(Z_s,p)}{\Sigma(Z_s)}\[\mu\(Z_s\)-r-\frac{\sigma(Z_s)b(Z_s){\rho}_WG_z(Z_s,s)}{G(Z_s,s)}\],\ t\leq s\leq T,
		\end{aligned}
	\end{equation} and the value function \begin{equation}
		V(x,y,z,t)=-x+G(z,t)y,
	\end{equation}
	satisfy Conditions of Theorem \ref{MMV verification} and are the solution to the MMV Problem \ref{MMV problem}, where $G(z,t)$ is the unique solution to the parabolic equation:
	$$\left\{\begin{aligned}
		&G_t+\alpha(z)G_z+\frac{1}{2}b(z)^2G_{zz}-\beta(z)b(z)^2\frac{G_z^2}{G}+\lambda(z)G=0,\\
		&G(z,T)=-1.
	\end{aligned}\right.$$
	Moreover, the optimal investment strategy $\pi^*$ can also be rewritten as
	$$\pi^*_s=-\frac{X^{\pi^*}_{s-}-x+2yG(z,t)}{\Sigma(Z_s)}\[\mu\(Z_s\)-r-\frac{\sigma(Z_s)b(Z_s){\rho}_WG_z(Z_s,s)}{G(Z_s,s)}\],\ t\leq s\leq T.$$
\end{theorem}

\proof{Proof.}
{ Using the fact that $G<0$ (we will prove it later in the proof of Theorem \ref{thm efficient frontier}), it is easy to see that the candidate control satisfies the HJBI Equation \eqref{MMV HJBI} and Conditions $(\romannumeral1)\sim(\romannumeral3)$ of Theorem \ref{MMV verification}.} Condition $(\romannumeral4)$ is also obvious. We only need to verify Condition $(\romannumeral5)$ and that $\pi^*,\eta_{1}^*,\eta_{2}^*$ and $\phi^*$ are admissible. 

{Under Conditions of Proposition \ref{pde condition}, we have $G=-e^{F}$ and $F\in C^{2+\beta,1+\beta/2}(H_T)$. Because $F$ and $F_z$ are bounded and $G_z=G\cdot F_z$, we have that $G$ and $G_z$ are bounded and $-G$ is bounded below from 0. Then, $\frac{G_z}{G}$ is bounded.} Using the Cauchy-Schwarz's inequality, we have
\begin{equation*}
	\begin{aligned} 
		&\frac{\big|\sigma(z)\(\mu(z)-r\)\big|}{\Sigma(z)}\leq \frac{\frac{1}{2}\sigma(z)^2+\frac{1}{2}\(\mu(z)-r\)^2}{\Sigma(z)}\leq \frac{1}{2}+\frac{1}{2}\lambda(z)\leq M,\\
		&\int_{\R}\[\frac{\gamma(z,p)\(\mu(z)-r\)}{\Sigma(z)}\]^2\nu\(\ud p\)\leq \frac{\frac{1}{2}\[\int_{\R}\gamma(z,p)^2\nu\(\ud p\)\]^2+\frac{1}{2}\(\mu(z)-r\)^4}{\Sigma(z)^2}\leq \frac{1}{2}+\frac{1}{2}\lambda(z)^2,\\
		&\int_{\R}\[\frac{\gamma(z,p)\sigma(z)}{\Sigma(z)}\]^2\nu\(\ud p\)\leq \frac{\frac{1}{2}\[\int_{\R}\gamma(z,p)^2\nu\(\ud p\)\]^2+\frac{1}{2}\sigma(z)^4}{\Sigma(z)^2}\leq \frac{1}{4}.
	\end{aligned}
\end{equation*}
Then, $\eta_{1}^*,\eta_{2}^*$ and $\phi^*$ (as functions of $Z_s$) satisfy Conditions of Lemma \ref{lem stochastic exponencial} and belong to $\mathcal{M}^t$. 

As for $\pi^*$, we define a jump diffusion process $R=\{R_s,s\in[t,T]\}$ and a function $\zeta$ as follows:
$$ R_s=\frac{X_s^{\pi^*}-x}{2yG(z,t)}+1,\ \zeta(z,s)=-\frac{\mu\(z\)-r-\frac{\sigma(z)b(z){\rho}_WG_z}{G}}{\Sigma(z)}. $$
As $\pi^*(s)=(X^{\pi^*}_{s-}-x+2yG(z,t))\zeta(Z_s,s)$, 
based on Theorem \ref{Y representation}, the jump diffusion process $R$ satisfies
\begin{equation}\label{dynamic of R}
	\ud R_s=R_{s}\zeta(Z_s,s)\(\mu\(Z_s\)-r\) \ud s+R_{s}\zeta(Z_s,s)\sigma\(Z_s\)\ud W^1_s+R_{s-}\zeta(Z_s,s)\int_{\mathbb{R}}\gamma\(Z_{s},p\) \tilde{N}\(\ud s,\ud p\). 
\end{equation}
{ As SDE \eqref{dynamic of R} is linear with bounded coefficients, using the useful moment estimate inequality for integrals with the L\'evy process (see, e.g., \citet{Protter2005}, Theorem 66 in Chapter 5), we have 
	$$\E_{x,y,z,t}\[\sup_{t\leq s\leq u}|R_s|^4\]\leq K\(1+\int_{0}^{u}\E_{x,y,z,t}\[\sup_{t\leq s\leq r}|R_s|^4\]\ud r\).$$
	The Gronwall's inequality yields $\E_{x,y,z,t}\[\sup\limits_{t\leq s\leq T}|R_s|^4\]<\infty,$ and then $$\E_{x,y,z,t}\[\sup_{t\leq s\leq T}|X_s^{\pi^*}|^4\]<\infty,$$
	As such,
	$${\E_{x,y,z,t}^{\eta,\phi}\[\sup_{t\leq s\leq T}|X_s^{\pi^*}|^2\]}\leq \(\E_{x,y,z,t}\[\(\frac{\ud \Q^{\eta,\phi}}{\ud \p}\)^2\]\E_{x,y,z,t}\[\sup_{t\leq s\leq T}|X_s^{\pi^*}|^4\]\)^{1/2}<\infty.$$}
Thus, the investment strategy $\pi^{\ast}$ is admissible. 

Finally, to show Condition $(\romannumeral5)$, by the definition of the admissible $\pi$, we have $$\E_{x,y,z,t}^{\eta,\phi}\[\sup_{t\leq s\leq T}|X_s^{\pi}|\]\leq \(\E_{x,y,z,t}^{\eta,\phi}\[\sup_{t\leq s\leq T}|X_s^{\pi}|^2\]\)^{1/2}<\infty.$$
On the other hand, using the Doob's maximal inequality for martingales, we obtain $$\begin{aligned}
	\E_{x,y,z,t}^{\eta,\phi}\[\sup_{t\leq s\leq T}|Y_s^{\eta,\phi}|\]\leq \(\E_{x,y,z,t}\[\(\frac{\ud \Q^{\eta,\phi}}{\ud \p}\)^2\]\E_{x,y,z,t}\[\sup_{t\leq s\leq T}|Y_s^{\eta,\phi}|^2\]\)^{1/2}
	\leq C\E_{x,y,z,t}\[|Y_T^{\eta,\phi}|^2\]<\infty.
\end{aligned}$$
As $G(z,t)$ is bounded, we see $$\E^{\eta,\phi}_{x,y,z,t}\[\sup_{t\leq s\leq T}|V(X_s,Y_s,Z_s,s))|\]\leq\E^{\eta,\phi}_{x,y,z,t}\[\sup_{t\leq s\leq T}|X_s^{\pi}|\]+C\E^{\eta,\phi}_{x,y,z,t}\[\sup_{t\leq s\leq T}|Y_s^{\eta,\phi}|\]<\infty,$$ which completes the proof.\hfill
\endproof 

\section{MV problem.}\label{section mv}

In this section, we consider the optimal investment problem under MV preferences in the same market model introduced in Section \ref{section model}. We use the embedding technique and the Lagrange multiplier method to transform the optimal investment problem into a series of auxiliary problems, which can be solved by the dynamic programming method. The dynamics of $X$ and $Z$ remain the same in Section \ref{section model}. Then, the dynamic of the forward wealth process is  
\begin{equation}\label{dynamic of Xt 2}
	\left\{\begin{aligned}
		&\ud X_t^{\pi}=\pi_{t}\(\mu\(Z_t\)-r\) \ud t+\pi_t\sigma\(Z_t\)\ud W^1_t+\pi_{t}\int_{\mathbb{R}}\gamma\(Z_{t},p\) \tilde{N}\(\ud t,\ud p\),\\ 
		&X_0=x.
	\end{aligned}\right.
\end{equation}
Again, we make the following assumption:
{\begin{definition}
		An investment strategy $\pi=\{\pi_s\}_{t\leq s\leq T}$ is admissible ($\pi\in \mathcal{\mathcal{B}}_{x,z,t}$) if and only if
		\begin{enumerate}[(\romannumeral1)]
			\item $\pi$ is predictable.
			\item $\pi$ is such that SDE (\ref{dynamic of Xt 2}) has a unique solution satisfying $$\E_{x,z,t}\[\sup_{t\leq s\leq T}|X_s^{\pi}|^2\]<\infty.$$
		\end{enumerate}
\end{definition}}
Let $U_\theta$ denote MV preferences:$$U_\theta(X)=\E_{x,z,t}\[X\]-\frac{\theta}{2}Var_{x,z,t}(X).$$
\begin{problem}\label{MV problem}
	Solve the stochastic control problem:
	\begin{equation}
		U(x,z,t)=\max_{\pi} U_\theta(X_T^{\pi}).
	\end{equation}
	If $\pi^*\in \mathcal{B}_{x,z,t}$ satisfies $U(x,z,t)=U_\theta(X_T^{\pi^*})$, then we call $\pi^*$ the optimal control of Problem \ref{MV problem} and $U(x,z,t)$ the value function of Problem \ref{MV problem}.
\end{problem}

{ Here, we take a useful method to simplify the problem (see, e.g., Chapter 6 in \citet{ZXY1999stochastic}). We embed the problem into tractable auxiliary problems.} First, we expose an extra condition $\E\[X_T^{\pi}\]=A$ for Problem \ref{MV problem} to get a subproblem:$$\max_{\pi}\  A-\frac{\theta}{2}\E_{x,z,t}\[\(X_T^{\pi}-A\)^2\],\ s.t.\ \E_{x,z,t}\[X_T^{\pi}\]=A.$$ Then, we vary $A$ and find the optimal $A$ among all the subproblems. For the subproblem, it is equivalent to find the optimal control of the problem:
\begin{equation}\label{MV subproblem}
	\min_{\pi}\  \E_{x,z,t}\[(X_T^{\pi}-A)^2\],\ s.t.\ \E_{x,z,t}\[X_T^{\pi}\]=A.
\end{equation}

Using the Lagrange method and setting a multiplier ($2c$), the problem \eqref{MV subproblem} can be transformed to be: For all $A$, find $\pi^*(c)$ to minimize $$\E_{x,z,t}\[\(X_T^{\pi}-(A+c)\)^2\]-2Ac-c^2,$$
and then find $c^*$ such that $\E_{x,z,t}\[X_T^{\pi^*(c^*)}\]=A$. After solving this optimization problem, we search for the optimal $A$ to maximize $U_\theta(X_T^{\pi^*(c^*)}(A))$.

Given $D\in \R$, we define an objective function and an auxiliary problem $$J^{\pi}_D(x,z,t)=\E_{x,z,t}\[(X_T^{\pi}-D)^2\].$$
\begin{problem}\label{MV auxiliary problem}
	Solve the stochastic control problem:
	\begin{equation}
		U_D(x,z,t)=\min_{\pi} J^{\pi}_D(x,z,t).
	\end{equation}
	If $\pi^*\in \mathcal{B}_{x,z,t}$ satisfies $U_D(x,z,t)=J^{\pi^*}_D(x,z,t)$, then we call $\pi^*$ the optimal control of Problem \ref{MV auxiliary problem} and $U_D(x,z,t)$ the value function of Problem \ref{MV auxiliary problem}.
\end{problem}

\subsection{Verification theorem and solution to Problem \ref{MV auxiliary problem}.}
We discuss the solution to the auxiliary Problem \ref{MV auxiliary problem} first. Similar to Section \ref{section MMV}, we use the dynamic programming method to solve it. The generator $\mathscr{B}$ of $X$ is:
\begin{equation}\label{MV generator}
	\begin{aligned}
		\mathscr{B}^{\pi}U(x,z,t)=&U_t+U_x\pi\[\mu\(z\)-r-\int_{\mathbb{R}} \gamma\(z,p\)\nu(\ud p) \]+U_za\(z\)\\
		&+\frac{1}{2}U_{xx}\pi^2\sigma(z)^2+\frac{1}{2}U_{zz}b(z)^2+U_{xz}\pi\sigma(z)\rho_Wb(z)\\
		&+\int_{\mathbb{R}}\[U(x+\pi\gamma(z,p),z,t)-U(x,z,t)\]\nu(\ud p).
	\end{aligned}
\end{equation}

We propose the corresponding verification theorem for Problem \ref{MV auxiliary problem} as follows:
\begin{theorem}\label{MV verification}
	Suppose that there exist a function $$U_D\in C^{2,2,1}\(\R\times\R\times[0,T)\)\cap C\(\R\times\R\times[0,T]\)$$ and a control $$\pi^*\in \mathcal{B}_{x,z,t}$$ such that 
	\begin{enumerate}[(i)]
		\item $\mathscr{B}^{\pi}U_D(x,z,t)\geq 0,\ \forall\(x,z,t\) \in \R \times \R \times [0,T),\ \forall\pi \in \R.$
		\item $\mathscr{B}^{\pi^*}U_D(x,z,t)= 0,\ \forall\(x,z,t\) \in \R \times \R \times [0,T).$
		\item $U_D(x,z,T)=(x-D)^2,\ \forall\(x,z\) \in \R \times \R.$
		\item $\E_{x,z,t}\[\sup_{t\leq s\leq T}|U_D(X_s^{\pi},Z_s,s))|\]<+\infty,\ \text{for} \ \forall\(x,z,t\) \in \R \times \R \times [0,T), \forall \pi \in \mathcal{B}_{x,z,t}$.
	\end{enumerate}
	Then, $\pi^*$ is the optimal control and $U_D(x,z,t)$ is the value function of Problem \ref{MV auxiliary problem}. 
\end{theorem}

\proof{Proof.}
The proof is similar to that of Theorem \ref{MMV verification}. For any $\pi\in \mathcal{B}_{x,z,t}$, we define a stopping time $\tau_N=T\wedge \inf\{t\geq 0: |X_t^{\pi}|\geq N,|Z_t|\geq N\}$. {Again, by Ito's formula (see, \citet{Protter2005}), we have \begin{equation}\label{dynkin formula 2}
		\E_{x,z,t}\[U_D(X_{\tau_N}^{\pi},Z_{\tau_N},\tau_N)\]=U_D(x,z,t)+\E_{x,z,t}\[\int_{t}^{\tau_N} \mathscr{B}^{\pi}U_D(X_{s-}^{\pi},Z_s,s)\ud s\].
\end{equation}}
Using Condition $(\romannumeral1)$, we have $$\E_{x,z,t}\[U_D(X_{\tau_N}^{\pi^*},Z_{\tau_N},\tau_N)\]\geq U_D(x,z,t).$$ 
Using Conditions $(\romannumeral3)$, $(\romannumeral4)$, and  Dominated convergence theorem, we obtain$$J_D^{\pi}(x,z,t)= \E_{x,y,z,t}\[\(X_T^{\pi^*}-D\)^2\]\geq U_D(x,z,t),$$ by letting $N\to \infty$. Thus, it indicates that
$$\min_{\pi} J_D^{\pi}(x,z,t)\geq U_D(x,z,t).$$
Using the same way, we have $$U_D(x,z,t)=J_D^{\pi^*}(x,z,t).$$ Thus,
$$U_D(x,z,t)=\min_{\pi} J_D^{\pi}(x,z,t)=J_D^{\pi^*}(x,z,t),$$
which completes the proof.\hfill
\endproof 

\begin{remark}
	Similar to Remark \ref{rmk MMV HJBI}, Conditions $(\romannumeral1)\sim(\romannumeral3)$ can be rewritten as the form:
	\begin{equation}\label{MV HJB}
		\left\{\begin{aligned}
			&\min_{\pi}\mathscr{B}^{\pi}U_D(x,z,t)=0,\\
			&U_D(x,z,T)=(x-D)^2,
		\end{aligned}\right.
	\end{equation}
	which is known as the Hamilton–Jacobi–Bellman (HJB) equation of Problem \ref{MV auxiliary problem} with the boundary condition. To obtain the explicit form of the solution to Equation \eqref{MV HJB}, we guess and verify the possible solution as what we have done in Section \ref{section MMV}. 
\end{remark}

Based on the boundary condition, we guess that $U_D$ has the form 
\begin{equation}
	U_D(x,z,t)=(x-D)^2H(z,t).
\end{equation} Substituting the special form of $U_D$ into the generator \eqref{MV generator}, we have
$$\begin{aligned}
	\mathscr{B}^{\pi}U_D(x,z,t)=&(x-D)^2H_t+2(x-D)H\pi\[\mu\(z\)-r-\int_{\mathbb{R}} \gamma\(z,p\)\nu(\ud p) \]+(x-D)^2H_za\(z\)\\
	&+H\pi^2\sigma(z)^2+\frac{1}{2}H_{zz}b(z)^2+2(x-D)H_z\pi\sigma(z)\rho_Wb(z)\\
	&+H\int_{\mathbb{R}}\[2\pi\gamma(z,p)(x-D)+\pi^2\gamma^2(z,p)\]\nu(\ud p),
\end{aligned}$$
which is a quadratic function of $\pi$. Again, we consider the first order condition and get the candidate optimal control
$$\pi^*(D)=-\frac{x-D}{\Sigma(z)}\[\mu\(z\)-r+\frac{\sigma(z)b(z){\rho}_WH_z}{H}\].$$
Substituting $\pi^*(D)$ into Equation \eqref{MV HJB}, the HJB equation becomes 
\begin{equation}\label{MV pde H}
	\left\{\begin{aligned}
		&H_t+\alpha(z)H_z+\frac{1}{2}b(z)^2H_{zz}-\(1-\beta(z)\)b(z)^2\frac{H_z^2}{H}-\lambda(z)H=0,\\
		&H(z,T)=1.
	\end{aligned}\right.
\end{equation}
By straightforward computation, we find that the solution to Equation \eqref{MMV pde G} and the solution to Equation \eqref{MV pde H} have the relationship $G=-\frac{1}{H}$.
Using the fact that $\frac{H_z}{H}=-\frac{G_z}{G}$, we observe that the candidate optimal control of the MV problem has a similar form as that of the MMV problem. 
\begin{theorem}\label{MV optimal control thm}
	Under Conditions of Proposition \ref{pde condition}, the optimal control of the auxiliary Problem \ref{MV auxiliary problem} is \begin{equation}\label{MV optimal control}
		\pi^*(D)(s)=-\frac{X_{s-}^{\pi^*(D)}-D}{\Sigma(Z_s)}\[\mu\(Z_s\)-r-\frac{\sigma(Z_s)b(Z_s){\rho}_WG_z(Z_s,t)}{G(Z_s,t)}\],\ t\leq s\leq T.
	\end{equation} The corresponding value function is $$U_D(x,z,t)=-(x-D)^2/G(z,t).$$
\end{theorem}

\proof{Proof.}
We just note that $H(z,t)=-\frac{1}{G(z,t)}$ is bounded, as $-G(z,t)$ is bounded below from 0. The remaining proof is exactly the same as that of Theorem \ref{MMV optimal control thm}. \hfill
\endproof 

\subsection{Solution to Problem \ref{MV problem}.}

Now, for any given $D\in\R$, we get the optimal control for the auxiliary Problem \ref{MV auxiliary problem}. For the original MV Problem \ref{MV problem}, we still need to find $c^*(A)$ such that $\E_{x,z,t}\[X_T^{\pi^*(A+c^*(A))}\]=A$ and $A^*$ that optimizes $U_{\theta}(X_T^{\pi^*(A+c^*(A))})$. Similar to \citet{trybula2019continuous}, we propose and prove the following theorem to give the solution. 
\begin{theorem}\label{MV optimal A}
	Recall that the jump diffusion process $R=\{R_s,s\in[t,T]\}$ defined in Theorem \ref{MMV optimal control thm} satisfies SDE:
	\begin{equation}\label{dynamic of Rt}
		\left\{\begin{array}{ll}
			&\ud R_s=R_{s}\zeta(Z_s,s)\(\mu\(Z_s\)-r\) \ud s+R_{s}\zeta(Z_s,s)\sigma\(Z_s\)\ud W^1_s\\
			&\quad\quad\quad+R_{s-}\zeta(Z_s,s)\int_{\mathbb{R}}\gamma\(Z_{s},p\) \tilde{N}\(\ud s,\ud p\),\\ 
			&R_t=1,\\
			&\zeta(z,s)=-\frac{\mu\(z\)-r+\frac{\sigma(z)b(z){\rho}_WH_z}{H}}{\Sigma(z)}.
		\end{array}\right.
	\end{equation}
	Then, we have $$c^*(A)=(A-x)\frac{\E_{z,t}[R_T]}{1-\E_{z,t}[R_T]},\  A^*+c^*(A^*)=x+\frac{1}{\theta}\cdot\frac{1-\E_{z,t}[R_T]}{Var_{z,t}[R_T]}.$$
\end{theorem}

\proof{Proof.}
As $$\pi^*(A+c)(s)=(X_{s-}^{\pi^*(A+c)}-(A+c))\zeta(Z_s,s),$$ it is easy to find that $$R_s=\frac{X_s^{\pi^*(A+c)}-(A+c)}{x-(A+c)}$$ satisfies SDE \eqref{dynamic of Rt} with an initial value equalling to one. Then, the Lagrange multiplier $c^*(A)$ can be solved by $$\E_{z,t}\[A+c^{\ast}+R_T\(x-(A+c^{\ast})\)\]=A.$$ Simple calculation yields  $c^*(A)=(A-x)\frac{\E_{z,t}[R_T]}{1-\E_{z,t}[R_T]}$. Now, we have the value function $$\begin{aligned}
	U_{\theta}(X_T^{\pi^*(A+c^*(A))})=&A-\frac{\theta}{2}\Big[\E_{z,t}\[\(A+c^*(A)+R_T\(x-(A+c^*(A))\)\)^2\]-A^2\Big]\\
	=&A-\frac{\theta}{2}\[(x-A)^2\frac{\E_{z,t}\[R_T^2\]-\E_{z,t}\[R_T\]^2}{(1-\E_{z,t}\[R_T\])^2}\].
\end{aligned}$$ It is a quadratic function of $A$, where the maximum point is $$A^*=x+\frac{1}{\theta}\cdot\frac{(1-\E_{z,t}\[R_T\])^2}{Var_{z,t}(R_T)}.$$Then, the optimal parameter for the auxiliary problem is $$D^*=A^*+c(A^*)=x+\frac{1}{\theta}\cdot\frac{1-\E_{z,t}\[R_T\]}{Var_{z,t}(R_T)}.$$\hfill
\endproof 

Based on Theorem \ref{MV optimal A}, the optimal control of Problem \ref{MV problem} is then given by 
\begin{equation}\label{MV optimal control final}
	\pi^*_s=-\frac{X_{s-}^{\pi^*}-x-\frac{1}{\theta}\cdot\frac{1-\E_{z,t}[R_T]}{Var_{z,t}[R_T]}}{\Sigma(Z_s)}\[\mu\(Z_s\)-r+\frac{\sigma(Z_s)b(Z_s){\rho}_WH_z(Z_s,s)}{H(Z_s,s)}\],\ t\leq s\leq T.
\end{equation}

\section{Main results.}\label{section relationship}

In Section \ref{section MMV} and Section \ref{section mv}, based on a jump diffusion financial market setting, we obtain the optimal investment strategies under MMV (with Conditions of Proposition \ref{pde condition} and Assumption \ref{assumption jumpsize}) and MV preferences separately. In this section, we compare the results of MMV and MV problems and propose that the optimal investment strategies under the two preferences are actually the same if and only if Conditions of Proposition \ref{pde condition} and Assumption \ref{assumption jumpsize} hold. Further, we illustrate some economic analysis about why or why not the two solutions may be different. 

\subsection{Connection between MMV and MV.}

We prove that under the setting in Theorem \ref{MV optimal A} and Assumption \ref{assumption jumpsize}, the optimal controls given by \eqref{MMV optimal control final} and \eqref{MV optimal control final} are the same. First, we give a crucial lemma on the analogy of the Feynman-Kac formula. The basic setting of the probability space remains the same as Section \ref{section model}.
{\begin{lemma}\label{Feynman-kac}
		Suppose that function $f(z,t)$ is of class $C^{2,1}(H_T)$ and satisfies  $$f_t+\frac{1}{2}b(z)^2f_{zz}+a(z)f_z+g(z)f=0.$$ Then, we have the representation $$\begin{aligned}
			f(z,t)=&\E_{z,t}\bigg[f(Z_T,T)\exp\bigg\{\int_{t}^{T}g(Z_s)\ud s\bigg\}K_t^T\bigg],
		\end{aligned}$$
		where $Z_s$ satisfies $$\ud Z_s=\[a(Z_s)-\rho_Wb(Z_s)\sigma_0(Z_s)\]\ud t+b(Z_s)\(\rho_W\ud W^1_t+\bar{\rho}_W\ud W^2_t\),$$
		$K_t$ satisfies $$\left\{\begin{array}{ll}
			&\ud K_t^s=K_t^s\sigma_0\(Z_s\)\ud W^1_s+K_t^{s-}\int_{\mathbb{R}}\gamma_0\(Z_{s},p\) \tilde{N}\(\ud s,\ud p\),\\ 
			&K_t^t=1,
		\end{array}\right.$$ and $g$, $\sigma_0$, $\gamma_0$ are bounded measurable functions.
	\end{lemma}
	
	\proof{Proof.}
	First,	 applying  Ito's formula to process $\{f(Z_u,u)\exp\{\int_{t}^{u}g(Z_u)\ud s\},u\in[t,T]\}$,we have $$\begin{aligned}
		&\ud \(f(Z_u,u)\exp\left\{\int_{t}^{u}g(Z_u)\ud s\right\}\)\\=&\exp\left\{\int_{t}^{u}g(Z_u)\ud s\right\}\bigg[\(f_t+\[a(Z_u)-\rho_Wb(Z_u)\sigma_0(Z_u)\]f_z+\frac{1}{2}b(Z_u)^2f_{zz}+g(Z_u)f\)du\\&+f_z(Z_u,u)b(Z_u)\(\rho_W\ud W^1_u+\bar{\rho}_W\ud W^2_u\)\bigg]\\
		=&\exp\left\{\int_{t}^{u}g(Z_u)\ud s\right\}\bigg[-\rho_Wb(Z_u)\sigma_0(Z_u)f_z(Z_u,u)\ud u+f_z(Z_u,u)b(Z_u)\(\rho_W\ud W^1_u+\bar{\rho}_W\ud W^2_u\)\bigg].
	\end{aligned}$$
	Then, 
	$$\begin{aligned}
		&\ud \(f(Z_u,u)\exp\left\{\int_{t}^{u}g(Z_u)\ud s\right\}K_t^u\)\\
		=&\exp\left\{\int_{t}^{u}g(Z_u)\ud s\right\}\bigg[K_t^uf_z(Z_u,u)b(Z_u)\(\rho_W\ud W^1_u+\bar{\rho}_W\ud W^2_u\)+f(Z_u,u)K_t^u\sigma_0\(Z_u\)\ud W^1_u\\
		&\quad\quad\quad+f(Z_u,u)K_t^{u-}\int_{\mathbb{R}}\gamma_0\(Z_{u},p\) \tilde{N}\(\ud u,\ud p\)\bigg].
	\end{aligned}$$ 
	Defining a stopping time $\tau_N=\inf\{s\geq t:|Z_s|\geq 0,|K_t^s|\geq N\}$
	and taking the expectation on both sides of the last equality, we obtain
	$$\begin{aligned}
		&\E_{z,t}\bigg[f(Z_{T\wedge\tau_N},T\wedge\tau_N)\exp\bigg\{\int_{t}^{T\wedge\tau_N}g(Z_s)\ud s \bigg\}K_t^{T\wedge\tau_N}\bigg] \\
		&=\E_{z,t}\[f(Z_t,t)\exp\left\{\int_{t}^{t}g(Z_s)\ud s\right\}K_t^t\] =f(z,t).
	\end{aligned}$$
	Because $f$ and $g$ are bounded, and $\E_{z,t}\[\sup_{t\leq s \leq T} \(K^s_t\)^2\]<+\infty$, using  Dominate convergence theorem, letting $ N\to \infty  $, we get $$f(z,t)=\E_{z,t}\bigg[f(Z_T,T)\exp\bigg\{\int_{t}^{T}g(Z_s)\ud s\bigg\}K_t^T\bigg].$$
	\hfill
	\endproof }

Now, we propose the main theorem that establishes the relationship between the solutions to MMV and MV problems. 
\begin{theorem}\label{optimal control concide}
	Under Conditions of Proposition \ref{pde condition} and Assumption \ref{assumption jumpsize}, the optimal investment strategies of MMV and MV problems are the same.
\end{theorem}

\proof{Proof.}
Recall that the two optimal controls under MMV and MV problems are $$-\frac{2Y_{s-}^{\eta^*,\phi^*}G(Z_s,s)}{\Sigma(Z_s)}\[\mu\(Z_s\)-r-\frac{\sigma(Z_s)b(Z_s){\rho}_WG_z(Z_s,s)}{G(Z_s,s)}\],$$ and $$-\frac{X_{s-}^{\pi^*}-x-\frac{1}{\theta}\cdot\frac{1-\E_{z,t}[R_T]}{Var_{z,t}[R_T]}}{\Sigma(Z_s)}\[\mu\(Z_s\)-r+\frac{\sigma(Z_s)b(Z_s){\rho}_WH_z(Z_s,s)}{H(Z_s,s)}\],\ t\leq s\leq T.$$
Based on Theorem \ref{Y representation}, we have $$2Y^{\eta^*,\phi^*}_sG(Z_s,s)=X^{\pi^*}_s-x+2yG(z,t),$$and then $$2Y^{\eta^*,\phi^*}_{s-}G(Z_s,s)=X^{\pi^*}_{s-}-x+2yG(z,t).$$
As $y=\frac{1}{2\theta}$, we only need to prove that
\begin{equation}
	\frac{1-\E_{z,t}[R_T]}{Var_{z,t}[R_T]}=-G(z,t).
\end{equation}
As we mentioned before, $G$ and $H$ satisfy $GH=-1$. Thus, it is sufficient to prove $$\E_{z,t}[R_T]=H(z,t),\ \E_{z,t}[R_T^2]=H(z,t).$$
{
	Let $$\tilde{R}_T=R_Te^{\int_{t}^{T}-\zeta(Z_s,s)\(\mu\(Z_s\)-r\) \ud s}.$$ Then, we see that $\tilde{R}_T$ satisfies SDE:
	$$\left\{\begin{array}{ll}
		&\ud \tilde{R_s}=\tilde{R}_{s}\zeta(Z_s,s)\sigma\(Z_s\)\ud W^1_s+\tilde{R}_{s-}\zeta(Z_s,s)\int_{\mathbb{R}}\gamma\(Z_{s},p\) \tilde{N}\(\ud s,\ud p\),\\ 
		&\tilde{R}_t=1.
	\end{array}\right.$$
	
	Now, we use Lemma \ref{Feynman-kac} to calculate $\E_{z,t}[R_T]$ and $\E_{z,t}[R_T^2]$. 
	
	First, simple calculation yields that $H$ satisfies the following equation: \begin{equation}\label{eq H 2}
		\begin{aligned}
			&H_t+H_z\(a+\rho_Wb\zeta\sigma\)+\frac{1}{2}b^2H_{zz}+H\zeta(\mu-r)\\=&H_t+\alpha(z)H_z+\frac{1}{2}b(z)^2H_{zz}-\(1-\beta(z)\)b(z)^2\frac{H_z^2}{H}-\lambda(z)H			=0,
		\end{aligned}
	\end{equation}
	where $\zeta,\mu,\sigma$ and $\gamma$ represent $\zeta(Z_s,s),\mu(Z_s),\sigma(Z_s)$ and $\gamma(Z_s,p)$, respectively. We use Lemma \ref{Feynman-kac} by letting $f=H$, $g=\zeta(\mu-r)$, and $K_t^s=\tilde{R_s}$. It is easy to check that, under Conditions of Proposition \ref{pde condition}, these functions and the process satisfy the Condition of Lemma \ref{Feynman-kac}. Then,  
	$$\begin{aligned}
		\E_{z,t}\[R_T\]=\E\bigg[\exp\bigg\{\int_{t}^{T}\zeta(\mu-r)\ud s \bigg\}\tilde{R_T}\bigg] = H(z,t).
	\end{aligned}$$
	Second, $H$ also satisfies the following equation:
	\begin{equation}\label{eq H 3}
		\begin{aligned}
			&H_t+H_z\(a+2\rho_Wb\zeta\sigma\)+\frac{1}{2}H_{zz}b^2+H\(\frac{\sigma^2b^2\rho_W^2}{\Sigma}\frac{H_z^2}{H^2}-\lambda\)\\=&H_t+\alpha(z)H_z+\frac{1}{2}b(z)^2H_{zz}-\(1-\beta(z)\)b(z)^2\frac{H_z^2}{H}-\lambda(z)H
			=0.
		\end{aligned}
	\end{equation}
	$R_T^2$ satisfies SDE:
	$$\left\{\begin{array}{ll}
		&\ud {R_s^2}={R}_{s}^2\[2\zeta(\mu-r)+\zeta^2\sigma^2+\zeta^2\int_{\mathbb{R}} \gamma^2\nu(\ud p)\]\ud s+{R}_{s}^22\zeta\sigma\ud W^1_s\\
		&\quad\quad\quad+{R}_{s-}^2\int_{\mathbb{R}}\(2\zeta\gamma+\zeta^2\gamma^2\) \tilde{N}\(\ud s,\ud p\),\\ 
		&\tilde{R}_t=1,
	\end{array}\right.$$
	where 
	$$\begin{aligned}
		2\zeta(\mu-r)+\zeta^2\sigma^2+\zeta^2\int_{\mathbb{R}} \gamma^2\nu(\ud p)
		=\frac{\sigma^2b^2\rho_W^2}{\Sigma}\frac{H_z^2}{H^2}-\lambda.
	\end{aligned}$$
	Using Lemma \ref{Feynman-kac} again for $f=H$, $g=\frac{\sigma^2b^2\rho_W^2}{\Sigma}\frac{H_z^2}{H^2}-\lambda$ and $K_t^s=R_s^2e^{\int_{t}^{s}- \(\frac{\sigma^2b^2\rho_W^2}{\Sigma}\frac{H_z^2}{H^2}-\lambda\)\ud u}$, we obtain
	$$\E_{z,t}\[R_T^2\]=\E\bigg[\exp\bigg\{\int_{t}^{T}\(\frac{\sigma^2b^2\rho_W^2}{\Sigma}\frac{H_z^2}{H^2}-\lambda\)\ud s\bigg\} K_t^T\bigg]=H(z,t).$$\hfill
	\endproof }

Until now, we prove that the optimal strategies of MMV and MV problems are the same under Conditions of Proposition \ref{pde condition} and Assumption \ref{assumption jumpsize}. However, to conclude that MMV and MV preferences perform exactly the same in the investment problem, we have to prove that they reach the same utility value under the optimal strategy. In the following lemma, we determine whether the investor's wealth process under the optimal strategy keeps in the domain of monotonicity of MV preferences (see Lemma 2.1 in \citet{maccheroni2009portfolio}) or not. This property is essential but somehow neglected by \citet{trybula2019continuous}.

{\begin{lemma}\label{MV in mono domain}
		Assume that we have Conditions of Proposition \ref{pde condition} for the existence of the solution to Equation \eqref{pde condition}. Under Assumption \ref{assumption jumpsize}, the wealth process under the optimal strategy \eqref{MV optimal control final} belongs to the domain of monotonicity of MV preferences. Conversely, if the wealth process under the optimal strategy \eqref{MV optimal control final} belongs to the domain of monotonicity of MV preferences for the initial value $(z,t)$, then Assumption \ref{assumption jumpsize} holds.
	\end{lemma}
	
	\proof{Proof.}
	For the first part of Lemma \ref{MV in mono domain}, we only need to prove $$X_T^{\pi^*(D^*)}\in \mathcal{G}_\theta:=\left\{f\in \mathcal{L}^2(\p):f-\E^{\p}\[f\]\leq \frac{1}{\theta} \right\}.$$
	By the definition of $R_t$, we have $$X_s^{\pi^*(D^*)}=x-\frac{1}{\theta}G(z,t)+\frac{1}{\theta}G(z,t)R_s.$$
	Letting $s=T$, we get $$X_T^{\pi^*(D^*)}=x-\frac{1}{\theta}G(z,t)+\frac{1}{\theta}G(z,t)R_T,$$ and
	$$\E_{x,z,t}\[X_T^{\pi^*(D^*)}\]=x-\frac{1}{\theta}G(z,t)+\frac{1}{\theta}G(z,t)H(z,t)=x-\frac{1}{\theta}G(z,t)-\frac{1}{\theta}.$$
	Then, $$X_T^{\pi^*(D^*)}-\E_{x,z,t}\[X_T^{\pi^*(D^*)}\]=-\frac{1}{\theta H(z,t)}R_T+\frac{1}{\theta}.$$
	If Assumption \ref{assumption jumpsize} holds, we have $R_T\geq0$ (see the representation of the stochastic exponential in Definition \ref{def exponential}), and then $X_T^{\pi^*(D^*)}-\E_{x,z,t}\[X_T^{\pi^*(D^*)}\]\leq\frac{1}{\theta}.$
	
	Now we consider the converse part of the Lemma \ref{MV in mono domain}. We could see that $X_T^{\pi^*(D^*)}\in \mathcal{G}_\theta$ is equivalent to $R_T\geq 0$ a.s.. Define $\tilde{R}_T=R_Te^{\int_{t}^{T}-\zeta(Z_s,s)\(\mu\(Z_s\)-r\) \ud s}.$ Then, we have $\p(R_T<0)=\p(\tilde{R}_T<0)$ and $\tilde{R}_T$ satisfies a linear SDE:
	$$\left\{\begin{array}{ll}
		&\ud \tilde{R_s}=\tilde{R}_{s}\zeta(Z_s,s)\sigma\(Z_s\)\ud W^1_s+\tilde{R}_{s-}\zeta(Z_s,s)\int_{\mathbb{R}}\gamma\(Z_{s},p\) \tilde{N}\(\ud s,\ud p\),\\ 
		&\tilde{R}_t=1.
	\end{array}\right.$$
	Similar to the proof of Lemma \ref{lem stochastic exponencial}, because $\zeta\sigma$ and $\zeta\gamma$ are bounded functions, the stochastic exponential $\tilde{R}$ satisfies the Condition in Lemma \ref{lemma memin} and then $\tilde{R}$ is a martingale. Thus, if $\tilde{R}_T\geq0$ a.s., then for any possible stopping time $\tau$, we have  $\tilde{R}_\tau\geq0$ a.s..
	
	If Assumption \ref{assumption jumpsize} does not hold, we let $\tau'=\inf\{s> t:\zeta(Z_s,s)\gamma(Z_s,\Delta L_s)\one_{\Delta L_{s}\neq 0}< -1\}$, where $\tau=\inf\{s> t:\zeta(Z_s,s)\gamma(Z_s,\Delta L_s)\one_{\Delta L_{s}\neq 0}= -1\}$. Then, by definitions of $\tau$ and $\tau'$, we have $\p(\tau'<\tau)>0$ and $\p(\tilde{R}_{\tau'}<0)>0$, which yields a contradiction.
	\hfill
	\endproof }
\vskip 5pt			
Finally, we discuss the relationship of the optimal controls and value functions in MMV and MV problems. Recall that the MMV preferences we propose in Section \ref{section model} is defined by $$V_\theta(X)=-\Lambda_\theta(X)-\frac{1}{2\theta},$$
where $$\Lambda_\theta(X)=\sup_{\Q\in \mathcal{Q}}\E^{\Q}\[-X-\frac{1}{2\theta}\frac{\ud \Q}{\ud \p}\].$$
However, we still need to prove that the optimal control and value remain the same when considering $\mathcal{Q}_0$ as the set of potential probability measures, as Remark \ref{MMV Q define} says. Let $$\Lambda_\theta^0(X):=\sup_{\Q\in \mathcal{Q}_0}\E^{\Q}\[-X-\frac{1}{2\theta}\frac{\ud \Q}{\ud \p}\]$$ and $V_\theta^0$ be the corresponding utility. We have the following theorem.
{\begin{theorem}
		With the initial value $(z,t)$, the optimal utilities of MMV and MV problems are the same if and only if Assumption \ref{assumption jumpsize} holds, even when we consider $\mathcal{Q}_0$ as the set of potential robust probability measures of MMV preferences. 
	\end{theorem}
	
	\proof{Proof.}
	We use $V_{MV},V_{MMV}$ and $V_{MMV}^0$ to represent the preferences when considering MV, MMV under $\mathcal{Q}^t$ and MMV under $\mathcal{Q}_0$. Then, with the initial value $(z,t)$, we have $V_{MMV}(X_T^{\pi})=-V(x,\frac{1}{2\theta},z,t)-\frac{1}{2\theta}$. 
	
	By the definition and properties of MMV preferences given in \citet{maccheroni2009portfolio}, for any admissible $\pi$, we have $$V_{MV}(X_T^{\pi})\leq V_{MMV}^0(X_T^{\pi})\leq V_{MMV}(X_T^{\pi}).$$
	If Assumption \ref{assumption jumpsize} holds, simple calculation yields $V_{MV}(X_T^{\pi^*})= V_{MMV}(X_T^{\pi^*})$, which means
	$$V_{MV}(X_T^{\pi^*})= V_{MMV}^0(X_T^{\pi^*})= V_{MMV}(X_T^{\pi^*}).$$ If Assumption \ref{assumption jumpsize} does not hold, by Lemma \ref{MV in mono domain}, we conclude that $X_T^{\pi^*}$ does not keep in the domain of monotonicity of MV preferences for some $t$ and then 
	\begin{equation}\label{eq ineq}
		\sup_{\pi}V_{MV}(X_T^{\pi})=V_{MV}(X_T^{\pi^*})<V^0_{MMV}(X_T^{\pi^*})\leq \sup_{\pi}V^0_{MMV}(X_T^{\pi})\leq \sup_{\pi}V_{MMV}(X_T^{\pi}),
	\end{equation}
	which completes the proof. 
	\vskip 5pt        
	It is worth mentioning that if Assumption \ref{assumption jumpsize} does not hold, we can only get the optimal strategy of MV problem but do not know the optimal strategy of MMV problem. However Inequality \eqref{eq ineq} directly shows that the optimal strategies and values of MMV and MV problems are different.
	\hfill
	\endproof }
\begin{remark}
	We cannot claim that the optimal values of MMV and MV problems are the same only due to the fact that $X_T^{\pi^*}$ falls in the domain of monotonicity (Lemma \ref{MV in mono domain}). On the one hand, we need to show that the choice of $\mathcal{Q}_0$ or $\mathcal{Q}^t$ does not affect the optimal strategy or the value function. More importantly, we have to prove that the optimal controls in MMV and MV problems are actually the same, which cannot be simply deduced by Lemma \ref{MV in mono domain}.
\end{remark}
{\begin{remark}
		There are simple examples under which Assumption \ref{assumption jumpsize} does not hold. For example, if we treat the stochastic factor $Z$ as a constant, then we have $$\zeta\gamma=-\frac{\mu-r}{\sigma^2+\gamma^2\nu(\R)}\gamma.$$
		By choosing $\mu$, $\sigma$, $\nu$ and $\gamma$ appropriately, $\zeta\gamma<-1$ can be achieved, which results in the violation of Assumption \ref{assumption jumpsize}. See Subsection \ref{subsec reasonableness} for more discussion on the reasonableness of Assumption \ref{assumption jumpsize}.
\end{remark}}

{In conclusion, we prove that in the optimal investment problem, MMV and MV preferences perform exactly the same if and only if Assumption \ref{assumption jumpsize} holds. The assumption is a necessary and sufficient condition. It also shows that the solutions of MMV and MV problems can be different in the financial market with jumps when Assumption \ref{assumption jumpsize} does not hold, which is a significantly new finding compared with \citet{trybula2019continuous} and \citet{SL2020note}. Thus, we conclude that the discontinuity of the financial market makes the difference between the investor's optimal strategies under MMV and MV preferences.}

\subsection{Efficient frontier.}

Following \citet{markowits1952portfolio}, \citet{LG2021RAIRO-OR} and etc., we derive the efficient frontiers (or efficient combinations in previous works) under MMV and MV preferences. Consider a bi-objective optimization problem:$$\max_{\pi\in \mathcal{B}_{x,z,t}} \(\E_{x,z,t}\[X_T^{\pi}\],-Var_{x,z,t}\(X_T^{\pi}\)\).$$
We cannot find a unique optimal point for this problem, but can find the equilibrium strategies, which are the solutions to $$\max_{\pi\in \mathcal{B}_{x,z,t}} \E_{x,z,t}\[X_T^{\pi}\]+\frac{\theta}{2}\[-Var_{x,z,t}\(X_T^{\pi}\)\],\ \forall \theta> 0.$$ Economically, $\theta$ represent the risk aversion level of the investor under MV preferences. Mathematically, $\frac{\theta}{2}$ represents the weight to determine the equilibrium. We call the following set the efficient frontier under MV preferences:
$$\mathcal{O}_{MV}=\left\{\(\E_{x,z,t}\[X_T^{\pi^*}\],Var_{x,z,t}\(X_T^{\pi^*}\)\)\bigg|\ \pi^* \  \text{is optimal strategy in Problem \ref{MV problem}},\ \forall \theta> 0 \right\}.$$
And the corresponding efficient frontier under MMV preferences is
$$\mathcal{O}_{MMV}=\left\{\(\E_{x,z,t}\[X_T^{\pi^*}\],Var_{x,z,t}\(X_T^{\pi^*}\)\)\bigg|\ \pi^* \  \text{is optimal strategy in Problem \ref{MMV problem}},\ \forall \theta> 0 \right\}.$$
The following theorem gives the explicit expression of both efficient frontiers.
\begin{theorem} \label{thm efficient frontier}
	Under Assumption \ref{assumption jumpsize}, the efficient frontier for both MMV and MV problems is 
	$$\left\{\Big(\E_{x,z,t}\[X_T^{\pi}\],Var_{x,z,t}\(X_T^{\pi}\)\Big)\bigg| \E_{x,z,t}\[X_T^{\pi}\]=x+\sqrt{\[-G(z,t)-1\]Var_{x,z,t}\(X_T^{\pi}\)} \right\}.$$
\end{theorem}

\proof{Proof.}
Fix $\theta>0$. By Theorem \ref{optimal control concide}, the efficient frontiers of MMV and MV problems are the same. By Theorem \ref{MV optimal A}, we use the relationship between $X_T^{\pi^*}$ and $R_T$ to calculate the efficient frontier. According to the proof of Lemma \ref{MV in mono domain}, we have \begin{equation}\label{eq efficient E}
	\E_{x,z,t}\[X_T^{\pi^*}\]=x-\frac{1}{\theta}G(z,t)+\frac{1}{\theta}G(z,t)H(z,t)=x-\frac{1}{\theta}G(z,t)-\frac{1}{\theta},
\end{equation}
and \begin{equation}\label{eq efficient Var}
	\begin{aligned}
		Var_{x,z,t}\(X_T^{\pi^*}\)=&\E_{x,z,t}\[\(\frac{1}{\theta}G(z,t)R_T+\frac{1}{\theta}\)^2\]\\
		=&\frac{1}{\theta^2}\[G(z,t)^2\E_{x,z,t}\[R_T^2\]+2G(z,t)\E_{x,z,t}\[R_T\]+1\]\\
		=&\frac{1}{\theta^2}\[G(z,t)^2H(z,t)+2G(z,t)H(z,t)+1\]\\
		=&-\frac{1}{\theta^2}\[G(z,t)+1\].
	\end{aligned}
\end{equation}
{The computation also indicates that $G(z,t)\leq -1$, because variance is always greater than or equal to 0. }Thus, combining \eqref{eq efficient E} and \eqref{eq efficient Var}, the efficient frontier is given by
$$\left\{\Big(\E_{x,z,t}\[X_T^{\pi}\],Var_{x,z,t}\(X_T^{\pi}\)\Big)\bigg| \E_{x,z,t}\[X_T^{\pi}\]=x+\sqrt{\[-G(z,t)-1\]Var_{x,z,t}\(X_T^{\pi}\)} \right\}.$$\hfill
\endproof 

\subsection{Economic impact of jump diffusion.}\label{subsec eco}

In this subsection, we analyze the economic impact of the jump diffusion part on the optimal strategy and value function. The basic idea is to find the change of the optimal strategy and value function when the jump risk is considered. We compare our results with \citet{trybula2019continuous}, where only Brownian motions are considered in the financial market. In general, our work is a generalization of \citet{trybula2019continuous} and we prove that MMV and MV preferences may be different in the jump diffusion market.

Recall that the optimal control in both MMV and MV problems is
$$\pi_s^*=-\frac{X_{s-}-x+\frac{1}{\theta}G(z,t)}{\Sigma(Z_s)}\[\mu\(Z_s\)-r-\frac{\sigma(Z_s)b(Z_s){\rho}_WG_z(Z_s,s)}{G(Z_s,s)}\],\ t\leq s\leq T,$$where $G(z,t)$ satisfies
\begin{equation}\label{MMV pde2}
	\left\{\begin{aligned}
		&G_t+\alpha(z)G_z+\frac{1}{2}b(z)^2G_{zz}-\beta(z)b(z)^2\frac{G_z^2}{G}+\lambda(z)G=0,\\
		&G(z,T)=-1,\\
		&\alpha(z)=a\(z\)-2\rho_W\frac{\mu\(z\)-r}{\Sigma(z)}\sigma(z)b(z),\\
		&\beta(z)=\frac{\bar{\rho}_W^2\sigma(z)^2+\int_{\R}{\gamma(z,p)^2\nu\(\ud p\)}}{\Sigma(z)}>0,\\
		&\lambda(z)=\frac{\(\mu\(z\)-r\)^2}{\Sigma(z)}>0.
	\end{aligned}\right.
\end{equation}
When the financial market is only driven by Brownian motions, i.e. $\gamma(z,p)\equiv 0$, it is easy to check that our results coincide with \citet{trybula2019continuous}. As Equation \eqref{MMV pde2} does not have an explicit solution, we simply assume that $\mu,\sigma$ and $\gamma$ are constants to show the change more visually. In this case, the function $G$ does not depend on the stochastic factor $Z$. Equation \eqref{MMV pde2} reduces to 
{$$\left\{\begin{aligned}
		&G_t+ \frac{\(\mu-r\)^2}{\sigma^2+\gamma^2\nu(\R)}G=0,\\
		&G(z,T)=-1,
	\end{aligned}\right.$$}
and the solution is $G(t)=-e^{\frac{\(\mu-r\)^2}{\sigma^2+\gamma^2\nu(\R)}(T-t)}$.
Then, the optimal strategy is
{$$\pi_s^*=\[x+\frac{1}{\theta}e^{\frac{\(\mu-r\)^2}{\sigma^2+\gamma^2\nu(\R)}(T-t)}-X_{s-}^{\pi^*}\]\frac{\mu-r}{\sigma^2+\gamma^2\nu(\R)},\ t\leq s\leq T.$$}
When market only contains Brownian motions, the optimal strategy is 
{$$\pi_s^*=\[x+\frac{1}{\theta}e^{\frac{\(\mu-r\)^2}{\sigma^2}(T-t)}-X_{s-}^{\pi^*}\]\frac{\mu-r}{\sigma^2},\ t\leq s\leq T.$$}

\begin{figure}[htbp]
	\centering
	\includegraphics[width=0.6\linewidth]{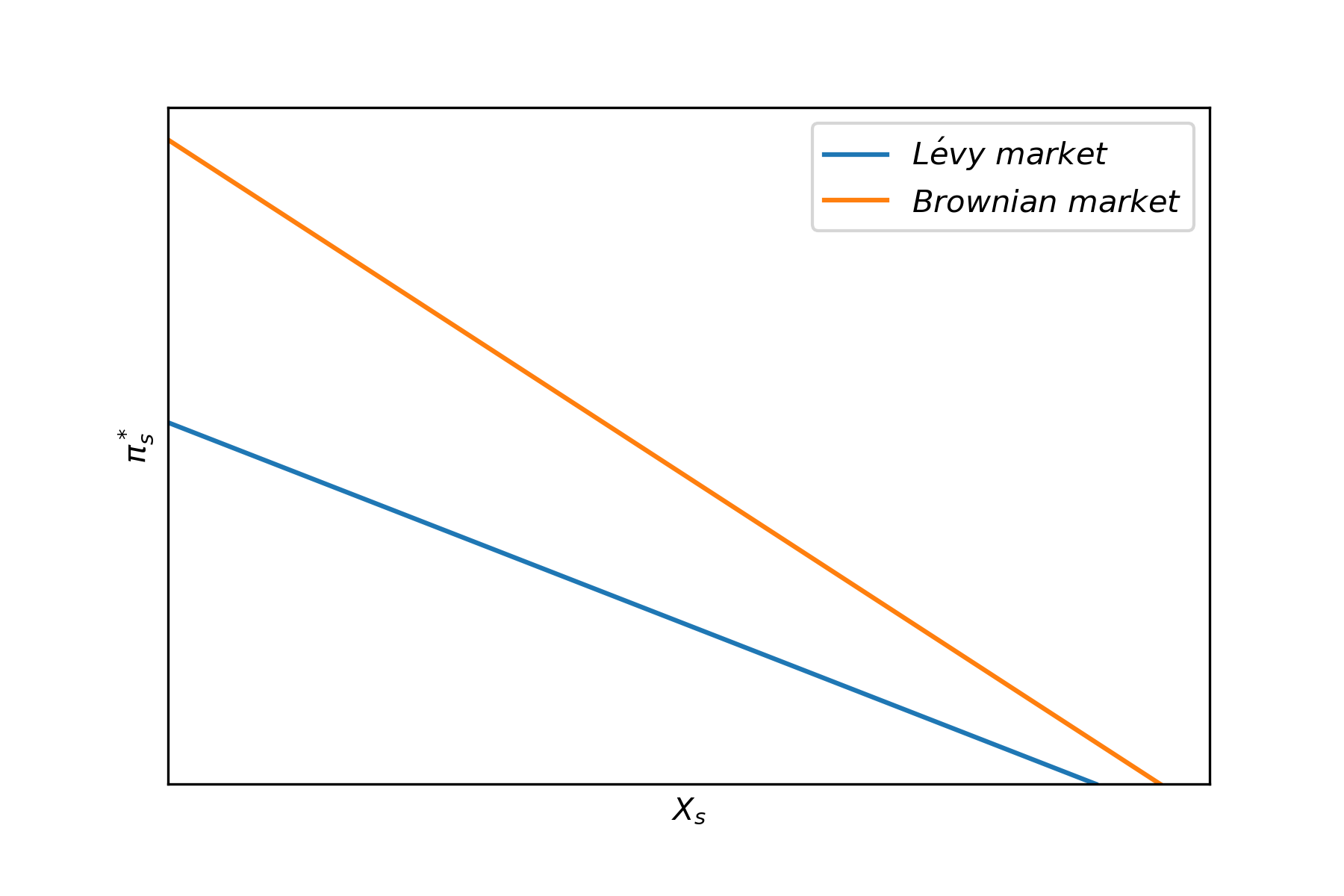}
	\caption{Optimal investment strategies in L\'{e}vy and Brownian financial markets}
	\label{picture_pi_X}
\end{figure}

Figure \ref{picture_pi_X} shows the change of the optimal strategy. When jump diffusion is considered, the investor tends to be more risk averse, and put less money in the risky asset. The phenomenon is reasonable because we can regard the Poisson jump part of the risky asset as an extra risk source. When the risk increases, the investor puts more money into the risk-free asset.

Moreover, we consider the change of the value function (or utility) when the market model changes. Based on the value function, the utility gained by the investor is
{$$\E_{x,z,t}\[X_T^{\pi^*}\]-\frac{\theta}{2}Var_{x,z,t}\(X_T^{\pi^*}\)=x-\frac{1}{2\theta}\[1-e^{\frac{\(\mu-r\)^2}{\sigma^2+\gamma^2\nu(\R)}(T-t)}\].$$}
Clearly, the existence of jump risk ($\gamma\neq 0$) leads to lower utility for the investor. It is reasonable as the jump part of the risky asset follows a martingale process with zero expectation, which does not provide any extra expected return but more extreme risks. When the variance increases and expected return declines, the investor's utility decreases. 

{\subsection{Reasonableness of Assumption \ref{assumption jumpsize}.}\label{subsec reasonableness}
	
	After discussing the equivalency of MMV and MV problems, an important question arises: When does Assumption \ref{assumption jumpsize} hold? In this subsection, we provide some numerical examples to show the reasonableness of Assumption \ref{assumption jumpsize}. Similar to Subsection \ref{subsec eco}, for the purpose of tractability we assume that $\mu$, $\sigma$ and $\gamma$ are constants, which means that the stochastic factor $Z$ is absent. Furthermore, we assume that the L\'evy process $L$ is a Poisson process $N$ with intensity $\lambda=\nu(\R)$.
	
	In this circumstance, the dynamic of risky asset $S$ is:$$\ud S_t=S_{t-}\(\mu\ud t+\sigma\ud W^1_t+\gamma\ud N_t\).$$
	\citet{honore1998pitfalls} has studied the problem of the parameter estimation in the above jump diffusion model. They focus on 18 liquid NYSE stocks and give the corresponding parameter estimations. We refer to their estimations to calculate the jump size we concern $$\zeta\gamma=-\frac{\mu-r}{\sigma^2+\gamma^2\nu(\R)}\gamma.$$
	By setting the risk-free rate $r = 0.04$, the estimation results are summarized in Table \ref{tab1}. 
	\begin{table}[htbp]\begin{center}\caption{Estimated parameters of different stocks.}\label{tab1}
			\begin{tabular}{lrrrrr}
				\hline  &   $\alpha$ &   $\sigma$ &  $\nu(\R)$ &     $\gamma$ &  $\zeta\gamma$ \\ 
				\hline   AXP     &  0.2887 &  0.3118 &  2.4754 & -0.068724 &   0.156934 \\
				CHV     &  0.2497 &  0.2463 &  2.6653 & -0.050006 &   0.155749 \\
				DD      &  0.1779 &  0.2348 &  1.2878 & -0.056067 &   0.130648 \\
				DOW     &  0.2044 &  0.2610 &  1.6259 & -0.061057 &   0.135311 \\
				EK      &  0.1254 &  0.2595 &  0.8759 & -0.082956 &   0.096561 \\
				GE      &  0.2389 &  0.2219 &  1.9257 & -0.052378 &   0.191077 \\
				GM      &  0.1460 &  0.2515 &  1.8357 & -0.053609 &   0.082924 \\
				IBM     &  0.1214 &  0.2351 &  0.7544 & -0.075313 &   0.102945 \\
				IP      &  0.1488 &  0.2598 &  0.4890 & -0.089990 &   0.137020 \\
				KO      &  0.2669 &  0.2392 &  1.6916 & -0.063026 &   0.223671 \\
				MMM     &  0.1549 &  0.2135 &  1.0352 & -0.063588 &   0.146807 \\
				MO      &  0.3067 &  0.2449 &  2.2174 & -0.057482 &   0.227783 \\
				MOB     &  0.2388 &  0.2458 &  2.0809 & -0.058235 &   0.171578 \\
				MRK     &  0.2656 &  0.2296 &  2.6653 & -0.042472 &   0.166567 \\
				PG      &  0.1782 &  0.2087 &  0.9792 & -0.058989 &   0.173588 \\
				S       &  0.1540 &  0.2656 &  1.2014 & -0.064337 &   0.097124 \\
				T       &  0.1676 &  0.1977 &  1.7907 & -0.049911 &   0.146251 \\
				XON     &  0.2135 &  0.1998 &  1.9356 & -0.047247 &   0.185291 \\
				DAX     &  0.2243 &  0.1377 &  4.6180 & -0.030621 &   0.242300 \\
				FTSE100 &  0.1522 &  0.1371 &  0.2963 & -0.088260 &   0.469223 \\
				SP100   &  0.2340 &  0.1396 &  2.1543 & -0.042089 &   0.350371 \\
				SP500   &  0.2464 &  0.1617 &  4.7382 & -0.039211 &   0.242077 \\
				 KFX     &  0.2359 &  0.1195 &  6.5222 & -0.024397 &   0.263150 \\
				\hline 
				\end{tabular}\end{center}
	\end{table}

	As we can see, $\zeta\gamma$ is positive in all the tested stocks because $\gamma$ is always negative. This implies that Assumption \ref{assumption jumpsize} holds almost in the real market. Furthermore, the absolute value $|\zeta\gamma|$ is significantly smaller than $1$, indicating that even when $\gamma>0$ or when $r$ fluctuates within a small range, Assumption \ref{assumption jumpsize} will remain valid.}

\section{Conclusion.}\label{section conclusion}

{ In this paper, we compare the optimal investment problems based on MMV and MV preferences under a jump diffusion and stochastic factor model. In the incomplete L\'{e}vy financial market, we prove that the optimal strategies and value functions of MMV and MV investment problems coincide if and only if a crucial assumption holds. Thus, the assumption is a sufficient and necessary condition for the consistency of MMV and MV preferences. When the key assumption violates, the solutions to the two investment problems can be different when there are jump risks in the financial market, which is an important addition to the literature. As \citet{SL2020note} prove that the optimal strategies of MMV and MV problems coincide in any continuous market, we provide an example to show that the discontinuity of asset prices can change such result. Thus, the difference between MMV and MV portfolio selections is due to the discontinuity of the market. Moreover, we provide some empirical evidences to illustrate the reasonableness of the assumption in real financial market. As recent developments in the MMV literature are also concerned with the trading constraints, we think that MMV portfolio selection problems when the market is discontinuous and the strategy is constrained can be further studied in the future.}

%
%
%

\section*{Acknowledgments.}
The authors acknowledge the support from the National Natural Science Foundation of China (Grant No.12271290, and No.11871036). The authors also thank the members of the group of Actuarial Sciences and Mathematical Finance at the Department of Mathematical Sciences, Tsinghua University for their feedback and useful conversations.

\bibliographystyle{plainnat}
\bibliography{jumpmmv-arxiv-revision.bib}
\end{document}